\documentclass[12pt]{article}
\usepackage{graphicx}
\usepackage{amssymb, cite}
\usepackage{setspace} 
\usepackage[hmargin=1.15in,vmargin=1.15in]{geometry} 

\newcommand{\be}{\begin{equation}}
\newcommand{\ee}{\end{equation}}
\newcommand{\bqn}{\begin{eqnarray}}

\newcommand{\eqn}{\end{eqnarray}}

\newcommand{\bd}{\begin{description}}
\newcommand{\ed}{\end{description}}

\newtheorem{stat}{}[section]

\def\bs{\begin{stat}}
\def\es{\end{stat}}
\def\ben{\begin{enumerate}}
\def\een{\end{enumerate}}

\def\bp{\noindent{\bf Proof}  \ \ \ }
\newcommand{\ep}{\hfill $\square$}

\makeatletter
\@addtoreset{equation}{section}

\makeatother

\begin{document}

\begin{center}
{\large {\bf PACKING  3-VERTEX PATHS  
IN CLAW-FREE GRAPHS
\\[2ex]
AND RELATED TOPICS}}
\\[4ex]
{\large {\bf Alexander Kelmans}}
\\[2ex]
{\bf University of Puerto Rico, San Juan, Puerto Rico, United States}
\\[0.5ex]
{\bf Rutgers University, New Brunswick, New Jersey, United States}
\\[2ex]
\end{center}

\begin{abstract}
A $\Lambda $-{\em factor} of a graph $G$ is a spanning 
subgraph of $G$ whose every component is a 3-vertex 
path. 
Let $v(G)$  be the number 
of vertices of $G$
 and $\gamma(G)$ the domination number of $G$.
A {\em claw} is a graph with four vertices and three edges incident to the same vertex. 
A graph is {\em claw-free} if it does not have an induced  subgraph isomorphic to a claw.
Our results include the following.
Let $G$ be a 3-connected claw-free graph, 
$x \in V(G)$, $e = xy \in E(G)$, and $L$ a 3-vertex path in $G$. Then 
$(a1)$ if $v(G) \equiv 0 \bmod 3$, then $G$ has  a 
$\Lambda $-factor  containing (avoiding) $e$,
$(a2)$ if $v(G) \equiv 1 \bmod 3$,  then $G - x$ has 
a $\Lambda $-factor, 
$(a3)$ if $v(G) \equiv 2 \bmod 3$, then $G - \{x,y\}$ has a $\Lambda $-factor,
$(a4)$ if $v(G) \equiv 0 \bmod 3$ and $G$ is either cubic or 4-connected, then $G - L$ has  a $\Lambda $-factor, 
$(a5)$ if $G$ is cubic with 
$v(G) \ge 6$ and 
$E$ is a set of three edges in $G$, 
then $ G - E$ has  a $\Lambda $-factor if and only if the subgraph induced by $E$ in $G$ is not a claw and not a triangle,
$(a6)$
if $v(G) \equiv 1 \bmod 3$,
then $G - \{v,e\}$ has a $\Lambda $-factor for every
vertex $v$ and every edge $e$ in $G$,
$(a7)$
if $v(G) \equiv 1 \bmod 3$, then there exist a 4-vertex path $\Pi$ and a claw $Y$  in $G$ such that $G - \Pi$ and $G - Y$ have $\Lambda $-factors, 
and
$(a8)$
$\gamma (G) \le \lceil v(G)/3 \rceil $ and if in addition 
$G$ is not a cycle and $v(G) \equiv 1 \bmod 3$, then $\gamma (G) \le \lfloor v(G)/3 \rfloor $. 
We also explore the relations between packing problems of a graph and its line graph to
obtain some results on different types of packings and 
 discuss  relations between 
$\Lambda $-packing and domination problems. 
\\[1ex]
  
\end{abstract}

\section{Introduction}
\label{Introduction}

\indent

We consider undirected graphs with no loops and 
no parallel edges unless stated explicitly. All notions and facts on graphs, that are  
used but not described here, can be found in \cite{BM,D,Wst}.

 Given a graph $G$ and a family ${\cal F}$ of non-isomorphic graphs,
an ${\cal F}$-{\em packing} of $G$ is a subgraph of $G$ 
whose every component is isomorphic to a member of 
${\cal F}$. An ${\cal F}$-packing $P$ of $G$ is called 
an ${\cal F}$-{\em factor} if $V(P) = V(G)$. 
The ${\cal F}$-{\em packing problem} is the problem of 
finding in $G$ an ${\cal F}$-packing having the maximum 
number of vertices.

If ${\cal F}$ consists of one graph $F$, then  
an ${\cal F}$-packing and an ${\cal F}$-factor are called simply an $F$-{\em packing} and an $F$-{\em factor}, respectively. Accordingly, the $F$-{\em packing problem} is the problem of finding in $G$ an $F$-packing having the maximum number of vertices or, equivalently, the maximum number of components.

If $F$ is a 2-vertex connected graph, then the 
$F$ packing problem is the classical matching problem  and a very beautiful and deep theory has been developed about this problem and its generalizations
(see, for example,   \cite{LP} as well as \cite{Kaneko}, 
\cite{Kstars}, and
\cite{LoPo}).
In particular, it is known that there is a polynomial-time algorithm for finding a maximum matching.
It turns out that  if $F$ is  a connected graph with at least three vertices, then  the $F$-packing problem is already $NP$-hard \cite{HK}.
Moreover, if $P_k$ is the $k$-vertex path, then for every 
$k \ge 3$ the $P_k$-packing problem turns out to be also  
 $NP$-hard for cubic graphs  \cite{K1}.

Let $\Lambda $ denote a 3-vertex path.
We will consider mainly the
$\Lambda $-packing problem.  This problem  is interesting for various reasons. Here are some of them.
\\[1ex]
\indent
$(R1)$ Path $\Lambda $ is  the smallest graph $F$, for which  the $F$-packing problem is $NP$-hard (even in the class of cubic graphs).
Although the $\Lambda $-packing problem  is $NP$-hard, i.e. possibly intractable in general,
it would be interesting to find some natural and non-trivial classes of graphs,  for which the problem is tractable, i.e. solvable in polynomial time
(e.g. 
{\bf \ref{r-regukargraphs}},
{\bf \ref{r-regukargraphsTreePack}}, 
{\bf \ref{2conclfr}}, and
{\bf \ref{Conclfr2endbplocks}} below). 
It is  also interesting to find polynomial-time algorithms  that provide a good approximation solution for the problem (e.g. 
{\bf \ref{km}} - {\bf \ref{r-regukargraphsTreePack}}
and 
{\bf \ref{eb(G)clfr}} below).
\\[1ex]
\indent
$(R2)$
Probably, one of the first non-trivial results in matching theory is Petersen's theorem (1891) stating that every cubic connected graph with at most two bridges has a perfect matching (see \cite{LP}). 
There are indications 
that a result of similar  nature may also be true for  the 
$\Lambda $-packing problem in the class of 3-connected graphs (see Problem {\bf \ref{Pr3con}} and theorems 
{\bf \ref{cubic3-conZTF}} - {\bf \ref{cubic3-conF}} below).
\\[1ex]
\indent
$(R3)$ 
It is known  \cite{Kaneko} that 
there is a  polynomial-time algorithm for the 
$\{P_3,P_4,P_5\}$-packing problem.
It can also be shown that 
a cubic 3-connected graph has a $\{P_3,P_4,P_5\}$-factor.
This fact for $\{P_3,P_4,P_5\}$-factors is analogous to 
Petersen's theorem for matchings mentioned above.
However, the complexity status of an $\{A,B\}$-packing 
problem for $A, B \in \{P_3,P_4,P_5\}$ and $A \ne B$ is not known. Some results in \cite{K3con-cub} (see also 
{\bf \ref{cubic3-conZTF}}  - {\bf \ref{cubic3-conF}}
below)) show that 
the $\Lambda $-packing problem for cubic 3-connected graphs is related to an $\{A,B\}$-packing  problem with
$A = P_3 = \Lambda $ and $B \in \{P_4,P_5\}$.
\\[1ex]
\indent
$(R4)$
The $\Lambda $-packing problem is also related to the
minimum domination problem in a graph (which is known to be $NP$-hard). Namely,  the size of a maximum $\Lambda $-packing in a graph $G$ can be used to give an upper bound for its domination number 
(see Section \ref{domination}). 
\\[1ex]
\indent
$(R5)$ The $\Lambda $-packing problem is also related to the
various problems on whether a graph $G$ has  a spanning subgraph $H$  of special type.  In the graph hamiltonicity theory $H$ is usually a Hamiltonian cycle or a Hamiltonian path. 
Obviously, the existence of  such a subgraph $H$ in a 
graph $G$ implies the existence  of  
a $\Lambda $-packing with  $\lfloor v(G)/3 \rfloor$ components.
For that reason,  various Hamiltonicity conjectures give rise to the corresponding $\Lambda $-factor problems or conjectures. (This was the original motivation to consider Problem {\bf \ref{Pr3con}} below.)
For example, in 1984 M. M. Mathews and D. P. Summer \cite {MS} conjectured that every 4-connected claw-free graph has a Hamiltonian cycle. 
Some results in the paper support this conjecture.
\\[1ex]
\indent
$(R6)$ Let $L(G)$ denote the line graph of $G$.
Then a vertex disjoint packing in $L(G)$ corresponds to an edge disjoint packing in $G$ and a vertex disjoint packing in $G$ corresponds to an induced vertex disjoint packing in $L(G)$. Since $L(G)$ is a claw-free graph, the study of 
the $\Lambda $-packing problem for claw-free graphs may allow to solve some problems on vertex and/or edge disjoint packings in graphs 
(see Section \ref{related-results}).
\\[1ex]
\indent
$(R7)$ The problem of packing induced 3-vertex paths in a claw-free graph is also related to the Hadwiger conjecture (see Section \ref{related-results}).
\\[1ex]
\indent
In Section \ref{notions} we give main notions and notation we use.
In Section \ref{preliminaries} we describe some known results and  open questions and outline main results  of the paper.
The formulations and proofs of the  main results on packings in claw-free graphs are given in Section \ref{clfree}. In Section \ref{domination} we discuss the relation between 
$\Lambda $-packing and domination problems and provide some bounds on the graph domination numbers based on some $\Lambda $-packing results.
Finally, in Section \ref{related-results}
we explore the relation between packing problems of a graph and its line graph to
obtain some results on different types of packings. 
We also discuss the induced $\Lambda $-packing problem and its relation with the Hadwiger conjecture.

\section{Main notions, notation, and simple observations}

\label{notions}

\indent

As usual, $V(G)$ and $E(G)$ denote the set of vertices and edges of $G$, respectively, and $v(G) = |V(G)|$, 
$e(G) = |E(G)|$. If $P$ is a path with the end-vertices $x$ and $y$, we put $End(P) = \{x,y\}$.
Given $X \subseteq E(G)$, let $\dot {X}$ denote the subgraph of $G$ induced by $X$.
Given $x \in V(G)$, let $N(x,G) = N(x)$ denote the set of vertices in $G$ adjacent to $x$. Let $Cmp(G)$ denote the set of components of $G$ and $cmp(G) = |Cmp(G)|$.
Let $\lambda (G)$ denote the maximum number of disjoint
3-vertex paths in $G$.
A vertex subset $X$ of $G$
is called a {\em domination set } in  $G$, if every vertex in $V(G) \setminus X$ is  adjacent to a vertex in $X$. 
Let $\gamma (G)$ denote the size of a minimum domination set  in  $G$; $\gamma (G)$ is called the {\em domination number} of $G$.
A {\em leaf} in a graph is a vertex of degree one. Let $Lv(G)$ denote the set of leaves in $G$ and $lv(G) = |Lv(G)|$.

A {\em claw} is a graph isomorphic to $K_{1,3}$, i.e. the graph with four vertices and three edges having a common end-vertex. A graph is called {\em claw-free} if it contains no induced claw.
A {\em net} is a graph obtained from a triangle with
three vertices $x_1$, $x_2$, and $x_3$ by adding three new vertices $z_1$, $z_2$, and $z_3$ and three new edges $x_1z_1$, $x_2z_2$, and $x_3z_3$.
A graph with one edge and two vertices is called a {\em match}.

A graph $G$ is {\em minimal 2-connected} if $G$ is 2-connected but $G - e$ is not 2-connected for every $e \in E(G)$. A {\em 2-frame} (or simply, a {\em frame}) of $G$ is a minimal 2-connected spanning subgraph of $G$.

Given a subgraph $S$ of a graph $G$, a vertex $x \in V(S)$ is a {\em boundary vertex} of $S$ if $x$  is adjacent to a vertex in $G - S $ and an {\em inner vertex} of $S$, otherwise.

A {\em block} of a connected graph $G$ is 
a maximal connected subgraph $H$ of $G$ such that $H - v$ is connected for every vertex $v$ of $H$, and so $H$ is either 2-connected or a match.
If $B$ has at most one boundary vertex,  then $B$ is called an {\em end-block} of $G$.
Let $eb(G)$ denote the number of end-blocks of connected graph $G$, and so if $eb(G) = 1$, then $G$ is either 
2-connected or a match.

We call a graph $H$ a {\em chain} if $H$ is connected and  has at most two end-blocks.
An {\em end-chain} of $G$ is a maximal  proper subgraph
$H$ of $G$ such that 
$H$ is a chain, 
every block of $H$  is a block of $G$ with at most two boundary vertices in $G$, and $H$ contains an end-block of $G$.
Obviously, a connected graph $G$ has an end-chain
if and only if $G$ has  at least three end-blocks.
 Also if $G$ has end-chains, then every end-block of $G$ is a subgraph (moreover, an end-block) of exactly one end-chain of $G$ and every end-chain of $G$ contains exactly one end-block of $G$.
 
 We call a graph $G$ a {\em $\Delta $-graph} if $G$ is  cubic and  every vertex of $G$ belongs to exactly one triangle, and so a $\Delta $-graph is a claw-free graph.
 
 We call a graph $G$  a {\em cactus} if $G$ is connected,  $G$ has at least three end-blocks,  and each end-chain 
of $G$ is a match.

Given a graph $G$, we write $G =  AxB$ if $A$  and $B$ are graphs, $V(A) \cap V(B) = \{x\}$, and 
$G  = A\cup B$, and so if $A$ and $B$ are connected graphs with at least two vertices, then $G$ is connected and $x$ is a cut-vertex of $G$.

We recall that a {\em $\Lambda $-packing} in a graph $G$ is a subgraph of $G$ whose every component is a 3-vertex graph and a {\em $\Lambda $-factor} in  $G$ is a spanning $\Lambda $-packing of $G$.
In addition, a  {\em $\Lambda $-packing} $P$ in  $G$ is called a {\em $\Lambda $-quasi-factor} of $G$ 
if $v(G) - v(P) \le 2$.
\\[1ex]
\indent 
 We will use the following simple facts. 
 \bs
\label{claw-freeAxB}
Let $G =  AxB$, where $A$ and $B$ are connected graphs with at least two vertices. Suppose that $G$ is claw-free.
Then the following holds.
\\[0.7ex]
$(a1)$ $N(x,A)$ and $N(x,B)$ induce complete subgraphs in $A$ and $B$, respectively.
\\[0.7ex]
$(a2)$ 
If $A$ is a block of $G$ and $v(A) \ge 3$, then 
$B - x$ is 
either 2-connected or a match,
and so $eb(B - x) = 1$.
\\[0.7ex]
$(a3)$ 
If $A$ is a block of $G$, $v(A) \ge 4$, and $xy \in E(A)$, then 
$A - \{x,y\}$ is a chain, and so $eb(A - x) \le 2$.
\\[0.7ex]
$(a4)$ 
If $A$ is a chain and $v(A) \ge 3$, then $A - x$ is also a  chain, and so $eb(A - x) \le 2$.
\\[0.7ex]
$(a5)$ 
If $A$ is an end-chain of $G$, $v(A) \ge 4$, and $xy \in E(A)$, then 
$A - \{x,y\}$ is connected and $eb(A - \{x,y\}) \le 3$.
\es

\section{Preliminaries and an outline of new results}

\label{preliminaries}
\indent

In \cite{K,KM} we gave an answer to 
the following natural question:
\\[1ex]
\indent
{\em How many disjoint 3-vertex paths must a cubic 
$n$-vertex graph have?}
\\[1ex]
\indent
Obviously, $\lambda (G) \le \lfloor v(G)/3 \rfloor $.
\bs 
\label{km}  
If $G$ is a cubic graph, then 
$\lambda (G) \ge \lceil v(G)/4  \rceil$ and at least  $v(G)/4$ disjoint 3-vertex paths in $G$ can be found in polynomial time.
\es

Obviously, if every component of $G$ is $K_4$, then
$\lambda (G) = v(G)/4$. 
Therefore the bound in {\bf \ref{km}} is sharp.
\\[.5ex]
\indent
Let ${\cal G}^3_2 $ denote the set of graphs with each vertex of degree  $2$ or $3$.
In \cite{K} we  gave (among other things) an answer
to the following  question:
\\[1ex]
\indent
{\em How many disjoint 3-vertex paths must an $n$-vertex  graph from ${\cal G}^3_2$ have?}

\bs 
\label{2,3-graphs} Suppose that $G \in {\cal G}^3_2$ and  
$G$ has no 5-vertex components.
Then $\lambda (G) \ge v(G)/4$ and at least 
$v(G)/4$ disjoint 3-vertex paths in $G$ can be found in polynomial time.
\es

From {\bf \ref{2,3-graphs}} it follows that 
every cubic graph $G$ has at least  $v(G)/4$ disjoint 3-vertex paths \cite{KM} because if $G$ is a cubic graph, then $G \in {\cal G}^3_2$ and $G$ has no 5-vertex components. 
\\[1ex]
\indent
In \cite{K} we  also gave a construction  that allowed  to prove the following:
\bs 
\label{extrgraphs1}
There are infinitely many 2-connected graphs in ${\cal G}^3_2$ $($and even subdivisions of cubic 3-connected graphs$)$ for  which the bound  in {\bf \ref{2,3-graphs}} is attained.
\es

Here are some packing results on regular graphs.
\bs {\em \cite{Kregular}}
\label{r-regukargraphs} 
Let $G$ be a $d$-regular graph with $d \ge 4$.
Then $\lambda (G) \ge v(G)/4$ and at least 
$v(G)/4$ disjoint 3-vertex paths in $G$ can be found in polynomial time.
\es

\bs {\em \cite{KMS}}
\label{r-regukargraphsTreePack} 
Let $T$ be a tree on $t$ vertices and let $\epsilon > 0$.
Suppose that $G$ is a $d$-regular graph and 
$d \ge \delta \ln \delta $, where $\delta = \frac{128t^3}{\epsilon ^2}$. Then $G$ contains at least 
$(1 - \epsilon) n/t $ vertex disjoint copies of $T$ and they can be found in polynomial time.
\es

There are infinitely many 2-connected cubic graphs $G$
having no $\Lambda$-quasi-factors.
Some of such graphs were constructed in \cite{Kcntrex} to provide 2-connected counterexamples to Reed's domination conjecture (see Section \ref{domination}).
In particular, a graph sequence $(R_k:  k \ge 3)$ in \cite{Kcntrex} is such that each $R_k$ is a cubic graph of connectivity two, $v(R_k) = 20k$, and
$\gamma (R_k)= (\frac{1}{3} + \frac{1}{60})v(R_k)$.
 Obviously, $\gamma (G) \le v(G) - 2 \lambda (G)$.
 Therefore 
$\lambda (R_k) \le  \frac{13}{40} v(R_k)$.
Questions arise whether there are 2-connected  cubic graphs with some additional properties and without  
$\Lambda$-quasi-factors.
For example,

\bs {\bf Problem.} 
\label{Pr2con,cub,bip,pl}  
Does every 2-connected, cubic, bipartite, and planar  graph have a $\Lambda$-quasi-factor?
\es

In \cite{K2con-cbp} we answered the question in 
{\bf \ref{Pr2con,cub,bip,pl}} by 
giving a construction that provides infinitely many  
 2-connected, cubic, bipartite, and planar  graphs 
without  $\Lambda$-quasi-factors. 
\\[1ex]
\indent
As to cubic 3-connected graphs, an old open question here is:

\bs {\bf Problem} {\em (A. Kelmans 1981).}
\label{Pr3con} 
Is the following claim true ?
\\[.5ex]
${\bf (P)}$ Every cubic 3-connected graph $G$ has 
a $\Lambda$-quasi-factor, i.e.
$\lambda (G) =   \lfloor v(G)/3 \rfloor $.
\es

In \cite{K3con-cub} we discuss Problem {\bf \ref{Pr3con}} and
show, in particular, that claim ${\bf (P)}$ in 
{\bf \ref{Pr3con}} is equivalent to some seemingly 
much stronger claims. Here are some results of this kind.

\bs {\em \cite{K3con-cub}}
\label{cubic3-conZTF}
The following are equivalent for cubic 3-connected graphs $G$$:$
\\[1ex]
$(z)$
$v(G) \equiv 0 \bmod 6$ $\Rightarrow$ $G$ has 
a $\Lambda $-factor,
\\[1ex]
 $(t)$ 
$v(G) \equiv 2 \bmod 6$ $\Rightarrow$ $G - \{x,y\}$ 
has a $\Lambda $-factor for some $x, y \in V(G)$, $x\ne y$, and
\\[1ex]
 $(f)$ 
$v(G) \equiv 4 \bmod 6$ $\Rightarrow$ $G - x$ 
has a $\Lambda $-factor for some $x \in V(G)$.
\es

\bs {\em \cite{K3con-cub}}
\label{cubic3-conZ}
The following are equivalent for cubic 3-connected graphs $G$ with $v(G) \equiv 0 \bmod 6$$:$
\\[1ex]
$(z0)$
$G$ has a $\Lambda $-factor,
\\[1ex] 
 $(z1)$ 
 for every $e \in E(G)$ there is a $\Lambda $-factor of 
 $G$ avoiding $e$, i.e. $G - e$ has a $\Lambda $-factor,
 \\[1ex]
 $(z2)$ 
 for every 
$e \in E(G)$ there is a $\Lambda $-factor of $G$ 
containing $e$,
 \\[1ex]
${(z3)}$
$G - X$ has 
a $\Lambda $-factor for every $X \subseteq E(G)$, 
$|X| = 2$,
and
\\[1ex]  
${(z4)}$
$G - L$ has a $\Lambda $-factor for every
3-vertex path $L$ in $G$.
\es

\bs {\em \cite{K3con-cub}}
\label{cubic3-conT}
The following are equivalent for cubic 3-connected graphs $G$ with $v(G) \equiv 2 \bmod 6$$:$
\\[1ex]
 $(t0)$ $G - \{x,y\}$ 
has a $\Lambda $-factor for some $x, y \in V(G)$, 
$x \ne y$, 
\\[1ex]
 $(t1)$ 
$G - \{x,y\}$ has a $\Lambda $-factor for some 
$xy \in E(G)$,
\\[1ex]
 $(t2)$ 
 $G - \{x,y\}$ 
has a $\Lambda $-factor for every $xy \in E(G)$,
and
\\[1ex]
$(t3)$
there  exists a 5-vertex path $W$ such that $G - W$ has 
a $\Lambda $-factor, and so $G$ has a $\{P_3,P_5\}$-factor.
\es

\bs {\em \cite{K3con-cub}}
\label{cubic3-conF}
The following are equivalent for cubic 3-connected graphs $G$ with $v(G) \equiv 4 \bmod 6$$:$
\\[1ex]
 $(f0)$
 $G - x$ 
has a $\Lambda $-factor for some 
$x \in V(G)$, 
\\[1ex]
 $(f1)$
 $G - x$ 
has a $\Lambda $-factor for every 
$x \in V(G)$, 
\\[1ex]
 $(f2)$ 
 $G - \{x, e\}$ 
has a $\Lambda $-factor for 
every $x \in V(G)$ and  every $e \in E(G)$,
and
\\[1ex]
$(f3)$ 
there exists a 4-vertex path $Z$ such that 
$G - Z$ has a $\Lambda $-factor, and so $G$ has 
a $\{P_3, P_4\}$-factor.
\es

There are some interesting results on the $\Lambda $-packing problem for  claw-free graphs.
Recall that a graph is called {\em claw-free} if it contains no induced subgraph isomorphic to a claw.
\bs {\em \cite{KKN}}
\label{2conclfr} 
Suppose that $G$ is a 2-connected claw-free graph.
Then
\\[0.7ex]
$(a1)$ if $v(G) \equiv 0 \bmod 3$, then $G$ has a $\Lambda $-factor,
\\[0.7ex]
$(a2)$ if $v(G) \equiv1 \bmod 3$, then $G - x$ has a $\Lambda $-factor for some $x \in V(G)$, and 
\\[0.7ex]
$(a3)$ if $v(G) \equiv 2 \bmod 3$, then $G -\{x,y\}$ has a 
$\Lambda $-factor for some $x,y \in V(G)$, $x \ne y$.
\\[0.7ex]
\indent
In every case a maximum $\Lambda $-packing can be found in polynomial time.
\es 

\bs {\em \cite{KKN}}
\label{eb(G)clfr}
Suppose that $G$ is a connected claw-free graph and  
$eb(G) \ge 2$. 
Then  $\lambda(G) \ge \lfloor (v(G) - eb(G) + 2)/3 \rfloor$,
this lower bound is sharp, and 
$\lfloor (v(G) - eb(G) + 2)/3 \rfloor $ disjoint 3-vertex paths in $G$ can be found in polynomial time.
\es

From {\bf \ref{eb(G)clfr}} we have, in particular:
\bs {\em \cite{KKN}}
\label{Conclfr2endbplocks} 
Suppose that $G$ is a  connected claw-free graph having exactly two end-blocks.
Then $\lambda(G) = \lfloor v(G)/3 \rfloor$ and 
a maximum $\Lambda $-packing can be found in polynomial time.
\es 

As we have mentioned in Section \ref{Introduction},
the $\Lambda $-packing problem remains $NP$-hard in the class of all graphs and even in the class of cubic graphs \cite{HK,K1}. It would be interesting to answer the following question.
\bs {\em {\bf Problem.}}
\label{claw-freeLpacNPhard}
Is the $\Lambda $-packing problem  $NP$-hard in the class of claw-free graphs ?
\es

In this paper (see Section \ref{clfree}) we give some more results on the $\Lambda $-packings in
claw-free graphs showing, in particular, to what extent the claims in 
{\bf \ref{cubic3-conZTF}} -  {\bf \ref{cubic3-conF}}
are true for claw-free graphs. Here are some of these results.
\\[1ex]
\indent
$(c1)$ 
If $G$ is a 2-connected  claw-free graph and
$v(G) \equiv 0 \bmod 3$, then for every edge  $e$ in $G$
there exists a $\Lambda $-factor of $G$ avoiding $e$, i.e. 
$G - e$ has a $\Lambda$-factor
(see {\bf \ref{clfree,2con-avoid-e}} and compare with {\bf \ref{cubic3-conZ}} $(z1)$ and {\bf \ref{2conclfr}} $(a1)$).
\\[1ex]
\indent
$(c2)$ 
If $G$ is a 3-connected  claw-free graph and
$v(G) \equiv 0 \bmod 3$, then for every edge $e$ in $G$ there exists a $\Lambda $-factor of $G$ containing $e$ 
(see {\bf \ref{clfree,3con,einL}} and compare with 
{\bf \ref{cubic3-conZ}} $(z2)$ and {\bf \ref{2conclfr}} $(a1)$).
\\[1ex]
\indent
$(c3)$
 If $G$ is a cubic 2-connected  claw-free graph with 
every vertex belonging to exactly one triangle
and 
$E$ is a set of two edges in $G$, 
then $ G - E$ has  a $\Lambda $-factor
(see {\bf \ref{clfree,2con-avoid(a,b)}} and compare with 
{\bf \ref{cubic3-conZ}} $(z3)$).
\\[1ex]
\indent
$(c4)$ 
If $G$ is a cubic 3-connected   claw-free graph with 
$v(G) \ge 6$ and 
$E$ is a set of three edges in $G$, 
then $ G - E$ has  a $\Lambda $-factor if and only if the subgraph induced by $E$ in $G$ is not a claw and not a triangle 
(see {\bf \ref{clfree,2con-avoid(a,b,c)}}).
\\[1ex]
\indent
$(c5)$ 
If $G$ is a cubic 3-connected  claw-free  graph 
or a 4-connected claw-free graph with 
$v(G) \equiv 0\bmod 3$, then for every 3-vertex path  $L$ in $G$ there exists a $\Lambda $-factor containing $L$, i.e. $G - L$ has a  $\Lambda $-factor
(see {\bf \ref{clfree,4con,deg3cntrL-}} and  {\bf \ref{2con,cubic,tr}}
and compare with {\bf \ref{cubic3-conZ}} $(z4)$).
\\[1ex]
\indent
$(c6)$ 
If $G$ is a 2-connected  claw-free graph and
$v(G) \equiv 2 \bmod 3$, then for every vertex $x$ in $G$ there exist  two edge $xy$ and $xz$ in $G$ such that 
$G - \{x,y\}$ and $G - \{x,z\}$ have 
$\Lambda$-factors
(see {\bf \ref{clfree,2con-avoid-e}} and  
{\bf \ref{clfree,2con,(x,y)-}}  and compare with 
{\bf \ref{cubic3-conT}} $(t1)$ and 
{\bf \ref{2conclfr}} $(a3)$).
\\[1ex]
\indent
$(c7)$ 
If $G$ is a 3-connected  claw-free graph and
$v(G) \equiv 2 \bmod 3$, then $G - \{x,y\}$ has 
a $\Lambda$-factor for every edge $xy$ in $G$
(see {\bf \ref{clfree,3con,(x,y)-}} and compare with 
{\bf \ref{cubic3-conT}} $(t2)$ and {\bf \ref{2conclfr}} $(a3)$).
\\[1ex]
\indent
$(c8)$
If $G$ is a 3-connected  claw-free 
graph and  $v(G) \equiv 2 \bmod 3$, then $G$ has a 5-vertex path $W$ such that $G - W$ has a $\Lambda $-factor, and so $G$ has 
a $\{P_3,P_5\}$-factor (see {\bf \ref{clfree,2con-avoid-e}} 
$(a2)$
and compare with {\bf \ref{cubic3-conT}} $(t3)$).
\\[1ex]
\indent
$(c9)$ 
If $G$ is a 2-connected  claw-free graph and
$v(G) \equiv 1 \bmod 3$, then $G - x$ has 
a $\Lambda$-factor for every vertex $x$ in $G$ 
(see {\bf \ref{clfree,2con,x-}} and compare with 
{\bf \ref{cubic3-conF}} $(f1)$ and {\bf \ref{2conclfr}} $(a2)$).
\\[1ex]
\indent
$(c10)$
If $G$ is a 3-connected  claw-free 
graph and  $v(G) \equiv 1 \bmod 3$,
then $G - \{x,e\}$ has a $\Lambda $-factor for every
vertex $x$ and every edge $e$ in $G$
(see {\bf \ref{clfree,3con,(x,e)-}} and compare with 
{\bf \ref{cubic3-conF}} $(f2)$).
\\[1ex]
\indent
$(c11)$
If $G$ is a 2-connected  claw-free 
graph and  $v(G) \equiv 1 \bmod 3$, then there exist a 4-vertex path $\Pi$ and a claw $Y$  in $G$ such that $G - \Pi$ and $G - Y$ have $\Lambda $-factors, and so $G$ has 
a $\{P_3,P_4\}$-factor and $\{P_3, Y\}$-factor
(see {\bf \ref{clfree,2con-avoid-e}} $(a2)$ and
{\bf \ref{G-Y}}
and  compare with {\bf \ref{cubic3-conF}} $(f3)$ and 
{\bf \ref{Ham}}).
\\[1ex]
\indent
$(c12)$
We show that the $\Lambda $-packing problem for a claw-free graph $G$ can be reduced in polynomial time to that for a special claw-free graph $K$ (called a cactus) with
$v(K) \le v(G)$ (see {\bf \ref{reduction}} and 
{\bf \ref{lambda(G)=}}).
\\[1ex]
\indent
$(c13)$ If $G$ is a 2-connected  claw-free graph,
then 
$\gamma (G) \le \lceil v(G)/3 \rceil $ and if in addition 
$G$ is not a cycle and $v(G) \equiv 1 \bmod 3$, then $\gamma (G) \le \lfloor v(G)/3 \rfloor $ (see {\bf \ref{gamma}}).

\section{Main results on claw-free graphs}
\label{clfree}

\indent

We will often use the following combination of 
{\bf \ref{2conclfr}} and {\bf \ref{Conclfr2endbplocks}}.
\bs {\em \cite{KKN}}
\label{eb<3}
Suppose that $G$ is a  connected claw-free graph  having at most two end-blocks.
Then $\lambda(G) = \lfloor v(G)/3 \rfloor$ and 
a maximum $\Lambda $-packing can be found in polynomial time.
\es 

Recall that  $G =  AxB$ is the union of graphs $A$ and $B$ with  $V(A) \cap V(B) = \{x\}$.
We need the following bounds on 
$\lambda  (AxB)$.
 \bs
\label{lambda(AxB)<}
 Let $G = AxB$, where $A$ and $B$ are connected graphs with at least two vertices. Then 
\\[0.5ex]
$(a0)$ if $v(A) \equiv 0\bmod 3$, then 
$\lambda (G) \le v(A)/3 + \lambda (B - x)$,
\\[0.5ex]
$(a1)$ if $v(A) \equiv 1\bmod 3$, then 
$\lambda (G) \le (v(A)-1)/3 + \lambda (B)$, and
\\[0.5ex]
$(a2)$ if $v(A) \equiv 2\bmod 3$, then 
$\lambda (G) \le (v(A) - 2)/3 + \lambda (B\cup xy)$, 
where $xy\in E(A)$.
 \es
 
 {\bf Proof.}
Let $S$ be a maximum $\Lambda $-packing in $G$.
\\[1ex]
\indent
${\bf (p1)}$ 
Suppose that $v(A) \equiv 0 \bmod 3$.
Let $S_1 = S \cap A$  and $S_2 = S - S_1$.
Then $\lambda (S_1) \le v(A)/3$.

Suppose that $x \in V(S_1)$. Then $\lambda (S_2) = \lambda (B - x)$, and so $\lambda (S) = \lambda (S_1) + \lambda (S_2) \le 
v(A)/3 + \lambda (B - x)$.

Now suppose that $x \not \in V(S_1)$. Then $\lambda (S_1) \le v(A)/3 - 1$. 
Let $S'_2 = S_2 \cap (B - x)$ and $S''_2 = S_2 - S'_2$.
Then 
$\lambda (S'_2) \le \lambda (B - x)$.
If $L$ is a 3-vertex path in $S''_2$, then $x \in V(L)$.  
Hence  $\lambda (S''_2) \le 1$. Therefore
$\lambda (S_2) = \lambda (S'_2) + \lambda (S''_2) \le 
\lambda (B - x) + 1$. Thus,
$\lambda (S) = \lambda (S_1) + \lambda (S_2) \le
(v(A)/3 - 1) + \lambda (B - x) + 1 = v(A)/3 + \lambda (B - x)$.
\\[1ex]
\indent
${\bf (p2)}$
Suppose that $v(A) \equiv 1\bmod 3$.
Let $S_1 = S \cap (A - x)$  and $S_2 = S - S'$.

Suppose that $\lambda (S_1) = v(A - x)/3$. 
Then  $\lambda (S_2) = \lambda (B)$, and so
$\lambda (S) = \lambda (S_1) + \lambda (S_2)  = 
v(A)/3 + \lambda (B - x)$.

Now suppose that $\lambda (S_1) \le v(A - x)/3 - 1$.
Let 
$S'_2 = S_2 \cap (B)$ and $S''_2 = S_2 - S'_2$.
Then 
$\lambda (S'_2) \le \lambda (B)$.
If $L$ is a 3-vertex path in $S''_2$, then $x \in V(L)$.  
Hence $\lambda (S''_2) \le 1$. Therefore
$\lambda (S_2) = \lambda (S'_2) + \lambda (S''_2) \le 
\lambda (B) + 1$. Thus,
$\lambda (S) = \lambda (S_1) + \lambda (S_2) \le
(v(A)/3 - 1) + \lambda (B) + 1 = v(A)/3 + \lambda (B)$.
\\[1ex]
\indent
${\bf (p3)}$ Finally, suppose that $v(A) \equiv 2\bmod 3$.
Let $S_1 = S \cap (A - \{x,y\})$ for some $xy \in E(A)$ and 
$S_2 = S - S_1$.
Let $S'_2 = S_2 \cap (B\cup xy)$ and $S''_2 = S_2 - S'_2$.
Then $\lambda (S'_2) \le \lambda (B\cup xy)$.
\\[0.7ex]
\indent
Suppose that $\lambda (S_1) = v(A - \{x,y\})/3$. 
Then $\lambda (S_2) = \lambda (B \cup xy)$, and so
$\lambda (S) = \lambda (S_1) + \lambda (S_2)  = 
v(A - \{x,y\})/3 + \lambda (B)$.
\\[0.7ex]
\indent
Suppose that $\lambda (S_1) = v(A - \{x,y\})/3 - 1$.
If $L$ is a 3-vertex path in $S''_2$, then $V(L) \subset 
V(A - S_1)$. Since $\lambda (S_1) = v(A - \{x,y\})/3 - 1$,
we have: $|V(A - S_1)| = 5$. Hence $\lambda (S''_2) \le 1$.
Therefore
$\lambda (S_2) = \lambda (S'_2) + \lambda (S''_2) \le 
\lambda (B) + 1$. Thus,
$\lambda (S) = \lambda (S_1) + \lambda (S_2) \le
(v(A - \{x,y\})/3 - 1) + \lambda (B\cup xy) + 1 =  
v(A - \{x,y\})/3 + \lambda (B\cup xy)$.
\\[0.7ex]
\indent
Now suppose that $\lambda (S_1) = v(A - \{x,y\})/3 - 2$.
If $L$ is a 3-vertex path in $S''_2$, then 
$V(L) \cap \{x,y\} \ne \emptyset $.
Hence $\lambda (S''_2) \le 2$.
Therefore
$\lambda (S_2) = \lambda (S'_2) + \lambda (S''_2) \le 
\lambda (B) + 2$. 
Thus, we have:
$\lambda (S) = \lambda (S_1) + \lambda (S_2) \le
(v(A - \{x,y\})/3 - 2) + \lambda (B) + 2 =  
v(A - \{x,y\})/3 + \lambda (B\cup xy)$.
\ep
\\[1ex]
\indent
It turns out that the end-chains of a claw-free graph have some special $\Lambda $-packing properties.

\bs
\label{chain2mod3}
Let $G$ be a connected claw-free graph, $C$ an end-chain of $G$, $v(C) \ge 3$, and $b$  the boundary vertex of $C$ {\em (and so $eb(C) \le 2$)}.
Then there exists an edge $bb'$ in $C$ such that
$\lambda (C - \{b,b'\}) = \lfloor v(C - \{b,b'\}/3) \rfloor $.
\es 

{\bf Proof} (uses {\bf \ref{eb<3}}).
Let $bx \in E(C)$. Since $G$ is claw-free,  $C - \{b,x\}$   is also claw-free. If there exists an edge $bb'$ in $C$ such that $eb(C - \{b,b'\}) \le 2$, then we are done by {\bf \ref{eb<3}}.

Let $B$ be the end-block of $C$ containing $b$, and so
$b$ is a boundary vertex of $B$. Since $G$ is claw-free,  $B$ and $B - \{b,x\}$ are also claw-free 
for $bx \in E(B)$ and $N(x,B)$ induces a complete subgraph $K$ in $B$.
\\[1ex]
${\bf (p1)}$
Suppose that $B$ has exactly one edge $bb'$. Since $v(C) \ge 3$, there is a (unique) block $B'$ in $C$ such that $b'$ is a boundary vertex of $B'$. 
If $e(B') = 1$, then 
$eb(C - \{b,b'\}) \le 2$  and we are done. If $e(B') \ge 2$, then $B'$ is 2-connected. Since $N(b',B')$ induces a complete subgraph in $B'$, clearly $B' - b'$ is either 2-connected or a match. 
Therefore again $eb(C - \{b,b'\}) \le 2$ and we are done.
\\[1ex]
${\bf (p2)}$
Now suppose
that $B$ has at least two edges, and so $B$ is 2-connected. 
Then $eb(B - \{b,x\}) \le 2$ for every edge $bx$ in $B$. 
Therefore if $B = C$, then we are done.
So we assume that $B \ne C$, and so $C$ has the end-block $D$ distinct from $B$.
If there is an edge $bb'$ in $B$ such that  
$eb(B - \{b,b'\})= 1$, then $eb(C - \{b,b'\}) \le 2$ and we are done.
So we assume that $eb(B - \{b,x\})= 2$ for every vertex 
$x$ in $K$. Then $v(K) \ge 2$. Let $B_1(x)$ and $B_2(x)$ be the two end-blocks of 
$B - \{b,x\}$, where $B_1(x)$ has no vertex adjacent to $b$.
Since $eb(B - \{b,x\})= 2$,  clearly $B_1(x)$ and $B_2(x)$ are 
2-connected for every edge $bx$ in $B$. 
Since $B \ne D$, clearly $eb(C - \{b,x\}) = 3$.
Since $G$ is claw-free, $C_x = C - (K - x)$ is also 
claw-free. 
Since $bx$ is an end-block of $C_x$, clearly 
$N(x, C_x - b)$ induces in $C_x - b$ a complete subgraph $K_x$. Now since  $B_1(x)$ has a vertex adjacent to $x$, clearly $K_x \subseteq B_1(x)$.
Since $B_2(x)$ is 2-connected, $K - x$ has an inner vertex  
$z$ of $B_2(x)$. Since $K - x$ is a complete graph and 
$B_2(x)$ is 2-connected, we have: 
$K - x \subseteq B_2(x)$.

First we assume that $v(K) \ge 3$, and so there is vertex 
$y$  in $K - \{x,z\}$. Then $y \in B_2(x)$, and so in 
$B -\{b,z\}$, vertex $x$ is adjacent to $B_1(x)$ and 
$B_2(x) - z$. Therefore $B - \{b,z\}$ is 2-connected, and so
$eb(B - \{b,z\})= 1$, a contradiction.

Now we assume that $v(K) = 2$, say, $V(K) = \{b_1 = x, b_2 = y\}$.
Let
$B_1$ be the subgraph of $B$ induced by $B_1(x) \cup x$
and 
$B_2 = B_2(x)$.
Then 
$N(b_1,C - b - b_2) = N(b_1,B_1)$,
$N(b_2,C - b - b_1)= N(b_2,B_2)$, and each $N(b_i, B_i)$ induces in $B_i$ a complete subgraph.
Let $B'_i = B_i - b_i$. 
Clearly, $B'_1 = B_1(b_1) = B_1(x)$.  
Since $B_1(x)$ is 2-connected, $B'_1$ is 2-connected.
Since $B_2 = B_2(x)$ is 2-connected, $B'_i$ is either 2-connected or a match.
Obviously, $B'_1$, $B'_2$, and $D$ are the three end-blocks of $C - \{b,b_1,b_2\}$. 
Let 
$C_i$ be the end-chain in $C - \{b,b_1,b_2\}$
containing  $B'_i$ and $c_i$ be the boundary vertex of  
$C_i$.
Let $C^i = C - \{b, b_i\}$ and $B'_i$ be the subgraph of $C$ induces by $C_i \cup b_i$.
\\[1ex]
${\bf (p2.1)}$
Suppose that $v(C_i) \equiv 0 \bmod 3$ for some $i \in \{1,2\}$, say, for $i = 1$.
Obviously, $C^1 - (C_1 - c_1)$  has exactly two end-blocks (namely, $B'_2$ and $D$).

Since $c_1$ is a vertex in $C^1 - (C_1 - c_1)$ and the boundary vertex of $C_i$, the neighborhood of 
$c_1$ in $C^1 - (C_1 - c_1)$ induces a complete subgraph in $C^1 - (C_1 - c_1)$.
Therefore $C^1 - C_1$ has also two end-blocks (namely, $B'_2$ and $D$).
By {\bf \ref{eb<3}}, 
$\lambda (C^1 - C_1) = \lfloor v(C^1 - C_1)/3 \rfloor $.
Since $C_1$ is a chain and $v(C_1) \equiv 0 \bmod 3$, 
by {\bf \ref{eb<3}}, 
$\lambda (C_1) = v(C_1)/3$.
Thus, 
$\lambda (C^1) = 
\lambda (C^1 - C_1) + v(C_1)/3 =
\lfloor v(C^1)/3 \rfloor $.
\\[1ex]
${\bf (p2.2)}$
Suppose that $v(C_1) \equiv 1 \bmod 3$ for some 
$i \in \{1,2\}$, say $i = 1$.
Since $c_1$ is the boundary vertex of end-chain $C_1$ in 
$C^1$, the neighborhood of 
$c_1$ in $C_1$ induces a complete subgraph in $C_1$. 
Therefore  $C_1 - c_1$ is a chain.
Since $v(C_1) \equiv 1 \bmod 3$, clearly 
$v(C_1 - c_1) \equiv 0 \bmod 3$. Therefore by {\bf \ref{eb<3}},
$\lambda (C_1 - c_1) = v(C_1 - c_1)/3$.
Obviously, $C^1 - (C_1 - c_1)$ is a chain.
By {\bf \ref{eb<3}},
$\lambda (C^1 - (C_1- c_1)) = 
 \lfloor v(C^1 - (C_1- c_1))/3 \rfloor $. 
Thus,
$\lambda (C^1) = 
\lambda (C^1 - C_1) + v(C_1- c_1)/3 =
\lfloor v(C^1)/3 \rfloor $.
\\[1ex]
${\bf (p2.3)}$
Finally, suppose that $v(C_1) \equiv 2 \bmod 3$ and 
$v(C_2) \equiv 2 \bmod 3$.
Let $C'_1$ denote the end-chain in $C^2$ containing 
$B'_1$. Then $C_1 \subset C'_1$ and 
$v(C'_1) = v(C_1) + 1 = 0 \bmod 3$. Now the arguments similar to those in ${\bf (p1)}$ shows that our claim is true.
\ep
\\[1ex]
\indent
Now we can improve bounds on $\lambda (G)$ in
{\bf \ref{lambda(AxB)<}}
when $G = AxB$ is claw-free and $A$ is an end-chain of $G$.
  \bs
\label{lambda(AxB)=}
 Let $G = AxB$, where $A$ and $B$ are connected graphs with at least two vertices. Suppose that $G$ is claw-free and $A$ is an end-chain of $G$. Then 
\\[0.5ex]
$(a0)$ if $v(A) \equiv 0\bmod 3$, then 
$\lambda (G)= v(A)/3 + \lambda (B - x)$,
\\[0.5ex]
$(a1)$ if $v(A) \equiv 1\bmod 3$, then 
$\lambda (G) = (v(A)-1)/3 + \lambda (B)$,
\\[0.5ex]
$(a2)$ if $v(A) \equiv 2\bmod 3$, then 
$\lambda (G) = (v(A)-2)/3 + \lambda (B\cup xy)$, 
where $xy$ is an edge in $A$.
 \es

{\bf Proof} (uses 
{\bf \ref{eb<3}},
{\bf \ref{lambda(AxB)<}},
and
{\bf \ref{chain2mod3}}).
Suppose that $v(A) \equiv 0\bmod 3$. Then by {\bf \ref{eb<3}},
$A$ has a $\Lambda $-factor $P$. Let $Q$ be a maximum
$\Lambda $-packing in $B - x$.
Then $P \cup Q$ is a $\Lambda $-packing in $G$ and
$\lambda (P \cup Q) = v(A)/3 + \lambda (B - x)$.
Therefore by {\bf \ref{lambda(AxB)<}} $(a0)$,
$\lambda (G)= v(A)/3 + \lambda (B - x)$.

Suppose that $v(A) \equiv 1\bmod 3$. Then by {\bf \ref{eb<3}},
$A-x$ has a $\Lambda $-factor $P$. Let $Q$ be a maximum $\Lambda $-packing in $B$.
Then $P \cup Q$ is a $\Lambda $-packing in $G$ and
$\lambda (P \cup Q) = (v(A)- 1)/3 + \lambda (B)$.
Therefore by {\bf \ref{lambda(AxB)<}} $(a1)$,
$\lambda (G)= (v(A)-1)/3 + \lambda (B)$.

Finally, suppose that  $v(A) \equiv 2\bmod 3$. Then by {\bf \ref{chain2mod3}}, there exists an edge $xy$ in $A$ such that
$A-\{x,y\}$ has a $\Lambda $-factor $P$. Let $Q$ be a maximum $\Lambda $-packing in $B$.
Then $P \cup Q$ is a $\Lambda $-packing in $G$ and
$\lambda (P \cup Q) = (v(A)- 2)/3 + \lambda (B\cup xy)$.
Therefore by {\bf \ref{lambda(AxB)<}} $(a2)$,
$\lambda (G)= (v(A)-2)/3 + \lambda (B\cup xy)$.
\ep
\\[1ex]
\indent
Theorem {\bf \ref{lambda(AxB)=}} suggests the following reduction procedure for claw-free graphs.
\\[1ex]
\indent
Let $G$ be a connected claw-free graph, $C$ an end-chain of $G$, and $c$  the boundary vertex of $C$.
Let us define a graph $\lfloor C \rfloor$ as follows: 
\\
if $v(C) \equiv 0 \bmod 3$ and $v(C) \ge 3$, then $\lfloor C \rfloor  = C$,
\\
if $v(C) \equiv 1 \bmod 3$ and $v(C) \ge 4$, then $\lfloor C \rfloor  = C - c$, and
\\
if $v(C) \equiv 2 \bmod 3$ and $v(C) \ge 5$, then $\lfloor C \rfloor  = C - \{c,c'\}$, where $cc'$ is an edge in $C$ such that $\lambda (C - \{c,c'\}) = \lfloor v(C - \{c,c'\}/3) \rfloor $
(see {\bf \ref{chain2mod3}}).

Obviously, $v(\lfloor C \rfloor)/3 =  \lfloor v(C)/3\rfloor $.
\\[1ex]
\indent
Recall that a graph  $G$ is called a {\em cactus} if $G$ is connected and has at least three end-chains and each end-chain has exactly two vertices.  
\\[1ex]
\indent
The following procedure for claw-free graphs allows to  either find a $\Lambda $-factor in a graph $G$ or to reduce 
the $\Lambda $-packing problem for $G$ to that for a cactus $K$ with $v(K) \le v(G)$.
\bs {\em {\bf Reduction.}}
\label{reduction}
Let $G$ be a connected claw-free graph.
\\[0.5ex]
$(s1)$ If $C_1$ is an end-chain of $G$ with $v(C_1) \ge 3$, then 
put $D_1 = \lfloor C _1 \rfloor$ and $G_1 = G - D_1$.
\\[0.5ex]
$(s2)$ We assume that $G_i$ and the sequence 
$(D_1, \dots , D_i)$ has already been defined for some $i \ge 1$. 

If $G_i$ has less than three end-chains or
every end-chain of $G_i$ has exactly two vertices, then stop and put $i = k$. 
Otherwise, let $C_{i+1}$ be an end-chain of $G$ with 
$v(C_{i+1}) \ge 3$.
Put $D_{i+1} = \lfloor C_{i+1} \rfloor $ and $G_{i+1} = G_i - D_{i+1}$.

The output of this procedure is 
$(D_1, \dots , D_k)$ and $G_k$.
\es

Obviously, Reduction {\bf \ref{reduction}} is a polynomial-time procedure.

Let $D^k = \cup \{D_i: i \in \{1, \ldots , k\}\}$.
Clearly all $D_i$'s are disjoint, and so
$\lambda (D^k) = \sum \{ \lambda (D_i): i \in \{1, \ldots , k\}\}$.

It follows that $G_k$ in Reduction {\bf \ref{reduction}} is either a claw-free chain or a claw-free cactus.
It is easy to show that if $G_k$ is a cactus, then $D^k$ and $G_k$ are uniquely defined; in this case let us denote 
$D_k$ by $D(G)$ and $G_k$ by 
$R(G)$.
\bs
\label{lambda(G)=}
Let $G$ be a connected claw-free graph and
$(D_1, \dots , D_k)$ and $G_k$ be the output of Reduction {\bf \ref{reduction}} applied to $G$. Let $Q$ be a maximum 
$\Lambda $-packing in $G_k$.
Then 
\\[0.5ex]
$(a1)$
each $D_i$  has a $\Lambda $-factor $P_i$, and so
$\lambda (P_i) = \lambda (D_i) = v(D_i)/3$,
\\[0.5ex]
$(a2)$
if $G_k$ is a chain, then
$P$ is a $\Lambda $-factor of $G$, 
\\[0.5ex]
$(a3)$
$P = Q \cup \{P_i: i \in \{1, \ldots , k\}\}$ is a maximum 
$\Lambda $-packing in $G$,
\\[0.5ex]
$(a4)$ if $G_k$ is not a chain, then
\\[0.3ex]
$\lambda (G) = \lambda (R(G)) + v(D(G)/3
\ge \lfloor (v(R(G)) - eb(R(G)) + 2)/3 \rfloor + 
v(D(G))/3 = l$,
\\
this lower bound is sharp, and $l$ disjoint  3-vertex paths in $G$ can be found in polynomial time.
\es

{\bf Proof} ~ 
(uses 
{\bf \ref{eb(G)clfr}}, 
{\bf \ref{eb<3}}, 
{\bf \ref{chain2mod3}}, and 
{\bf \ref{reduction}}).
We prove $(a1)$ and $(a2)$. 
By Reduction {\bf \ref{reduction}}, each $D_i$ has at most two end-blocks and $v(D_i) \equiv 0 \bmod 3$.
Since $G$ is claw-free, each $D_i$ is also claw-free.
By {\bf \ref{eb<3}}, each $D_i$ has a 
$\Lambda $-factor $P_i$, and so 
$\lambda (P_i) = \lambda (D_i) = v(D_i)/3$. Therefore $(a1)$ holds. 
If $G_k$ is a chain, then by the same reason, $Q$ is a 
$\Lambda $-factor of $G_k$. Then  
$P = Q \cup \{P_i: i \in \{1, \ldots , k\}\}$ is a  
$\Lambda $-factor of $G$, and so $(a2)$ holds. 
Now $(a3)$ follows from {\bf \ref{lambda(AxB)=}} and 
$(a4)$ follows from $(a3)$ and 
{\bf \ref{eb(G)clfr}}.
\ep
\\[1ex]
\indent
From {\bf \ref{reduction}} and {\bf \ref{lambda(G)=}} it follows that Problem {\bf \ref{claw-freeLpacNPhard}} is equivalent to
\bs {\em {\bf Problem.}}
Is $\Lambda $-packing problem $NP$-hard for claw-free cacti ?
\es

Now we describe an infinite class of sub-cubic claw-free graphs with no $\Lambda $ factors. This class includes infinitely many cacti. We will use this description to establish some 
$\Lambda $-packing properties of $\Delta $-graphs (see 
{\bf \ref{clfree,2con-avoid(a,b,c)}}).
\\[0.5ex]
\indent
Let ${\cal S}$ denote the set of graphs $S$ with the following properties:
\\[.5ex]
$(\alpha 1)$ $S$ is connected,
\\[.5ex]
$(\alpha 2)$ every vertex in $S$ has degree at most 3,
\\[.5ex]
$(\alpha 3)$ every vertex in $S$ of degree 2 or 3 belongs to exactly one triangle, 
and
\\[.5ex]
$(\alpha 4)$ $S$ has at least  three leaves.
\bs
\label{A}
 If $S \in {\cal S}$, then $S$ has no $\Lambda $-factor.
 \es

\bp Let $S \in {\cal S}$. If $v(S) \not \equiv 0 \bmod 3$, then our claim is obviously true. So we assume that  $v(S) \equiv 0\bmod 3$. By $(\alpha 3)$,   $v(S) \equiv lv(S) \bmod 3$, and so $lv(S) \equiv 0\bmod 3$. Obviously, it is sufficient to prove our claim for 
$S\in {\cal S}$ with property
$(\alpha '4)$: $lv(A) = 3$.
We prove our claim by induction on $v(G)$. The smallest graph in ${\cal S}$ is a net $N$ with $v(N) = 6$ and our claim is obviously true for $N$. So let $v(S) \ge 9$.
Suppose, on the contrary, that $S$ has a $\Lambda $-factor $P$. Let $v$ be a leaf of
$S$ and $vx$ the edge incident to $v$. 
Since $P$ is a  $\Lambda $-factor in $S$, it has a component $L = vxy$, and so $P - L$ is a $\Lambda $-factor in $S - L$ and $d(x, S) \ge 2$.
By property $(c3)$, $x$ belongs to a unique triangle $xyz$ in $A$ and $d(x,a) = 3$.
If $d(z,S) = 2$, then $z$ is an isolated vertex in $S - L$, and so $P$ is not a  $\Lambda $-factor in $S$, a contradiction.
Therefore by $(c2)$, $d(z,S) = 3$. Hence $z$ is a leaf in 
$S - L$, and so $lv(S - L) = 3$.
Therefore $S - L$ satisfies $(\alpha 2)$, $(\alpha  3)$, and 
$(\alpha '4)$.

Suppose that $G - L$ is not connected and that the three leaves 
do not belong to a common component. Then $S - L$ has a component $C$ with $v(C) \not \equiv 0 \bmod 3$, and so $S - L$ has no $\Lambda $-factor, a contradiction.

Finally, suppose that $S - L$ has a component $C$ containing all three leaves of $S - L$. Then $C \in {\cal S}$ and 
$v(C) < v(S)$. By the induction hypothesis, $C$ has no $\Lambda $-factor. Therefore $S - L$ also has no $\Lambda $-factor, a contradiction.
\ep
\\[1ex]
\indent
Recall that  a  {\em frame} of $G$ is a minimal 
2-connected spanning subgraph of $G$.
\\[1ex]
\indent
We need the following procedure from \cite{KKN} that
provides a frame of a 2-connected graph. 
This procedure was used in \cite{KKN} to prove {\bf \ref{2conclfr}}.
\bs {\em {\bf Procedure ${\cal E}$.}}
\label{ProcedureE}
Let $G$ be a 2-connected graph. We define
 sequences ${\cal A} = (A_0, \ldots , A_r)$ and  ${\cal G} =(G_0, \ldots , G_r)$ recursively, where each $A_i$ and each $G_i$ is a subgraph of $G$$:$
\\[1ex]
$(s1)$ Let $A_0$ be a longest cycle in $G$ and 
$G_0 = A_0$. 
\\[1ex]
$(s2)$ Assuming that the sequences $(A_0, \ldots , A_{i-1})$ and  $(G_0, \ldots , G_{i-1})$ are already defined,
let $A_i$ be a longest path in $G$ with the property
 
${\bf (E_i)}$$:$
$e(A_i) \ge 2$ and $G_{i-1} \cap A_i = End(A_i)$. 
\\
Put $G_i = G_{i-1} \cup A_i$. 
\\[1ex]
$(s3)$ Let $r$ be the minimum positive integer such that $G$ has no path $A_{r+1}$ with property ${\bf (E_{r+1})}$.
\es

If $G$ is a 2-connected graph, then we put
$F(G) = G_r$, $A(G) = A_r$, and ${\cal A}(G) = {\cal A}$
in  Procedure ${\cal E}$.
Clearly, every 2-connected graph has a frame.
\\[1ex]
\indent
It is easy to see the following.
\bs 
\label{frame}
Let $G$ be a 2-connected graph.
Then 
$F(G)$ is a frame of  $G$ and $F(G)$ is a Hamiltonian cycle of $G$ if and only if $r = 0$.
\es

We will also need the following modification of  
Procedure ${\cal E}$. 
Recall that an edge $e = xy$ is a {\em chord} of a cycle $C$ in $G$ if $e \not \in E(C)$ and $x,y \in V(C)$.
\bs {\em {\bf Procedure $\cal E$'.}}
\label{ProcedureE}
Let $G$ be a 2-connected graph and  $e \in E(G)$. 
Let Procedure ${\cal E}'$ be obtained from Procedure ${\cal E}$ by replacing the first step 
\\[0.7ex]
$(s1)$  
Let $A'_0 = G'_0$ be  a longest cycle in $G$.

by
\\[0.7ex]
$(s'1)$ Let $A'_0$ be a longest cycle  among all cycles $C$ in $G$ such that  edge $e$ is either in  $C$ or  is a chord of $C$ and let $G'_0 = A'_0$. 
\es

Since $G$ is 2-connected, $G$ has a cycle containing 
$e$. Therefore a cycle $A'_0$  exists.
\\[1ex]
\indent
It turns out \cite{KKN} that applied to a 2-connected claw-free graph $G$,  Procedure ${\cal E}$ provides a frame
$F(G) = G_r$  of $G$ and its ear-assembly with very useful properties.  

Recall that a {\em claw-free frame} of $G$ is a minimal 
2-connected claw-free spanning subgraph of $G$.
Clearly, every 2-connected claw-free graph has a claw-free frame.

\bs {\em \cite{KKN}}
\label{frame,claw-free}
Let $G$ be a  2-connected claw-free graph 
and $G$ not a cycle.
Let $F = F(G)$ and $A = A(G)$ from Procedure ${\cal E}$.
Then
\\[0.5ex]
$(f1)$ $F$ is a frame of $G$ with the maximum vertex degree three,
\\[0.5ex]
$(f2)$ $G$ has a unique matching $M$ such that
$F_c = F \cup M$ is a claw-free frame of $G$ with the maximum vertex degree three, and so
every vertex of degree three belongs to a unique triangle in $F_c$ and every vertex of every triangle in $F_c$ has degree three in $F_c$, 
\\[0.5ex]
$(f3)$ $F_c - A$ is a claw-free frame of $G - A$
$($and so $G - A$ is 2-connected and claw-free$)$
{\em (put $F_c = F_c(G)$)}, 
\\[0.5ex]
$(f4)$
if $P$ is a maximum $\Lambda $-packing of $A$,
then $G - P$ is a 2-connected claw-free graph, and
\\[0.5ex]
$(f5)$ the above claims are also true for Procedure 
${\cal E}'$.
\es
 
Obviously,  even the first steps in Procedures ${\cal E}$ and ${\cal E}'$  
are $NP$-hard. 
However, there are modifications  of these procedures
which find 
$F(G)$, $F_c(G)$, $A(G)$, and ${\cal A}(G)$ with properties  in {\bf \ref{frame,claw-free}} 
in polynomial time for every 2-connected claw-free graph $G$.

\bs 
\label{clfree,2con-avoid-e}
Suppose that $G$ is a 2-connected  claw-free graph.
\\[1ex]
$(a1)$ If 
 $v(G) \equiv 0 \bmod 3$ and $e \in E(G)$,
then  $G - e$ has a $\Lambda $-factor.
\\[1ex]
$(a2)$ If 
 $v(G) \equiv k \bmod 3$, where $k \in \{1,2\}$,
then  $G$ has a $k$-vertex path $P_k$ and a $(k+3)$-vertex path 
$P_{k+3}$ such that
$G - P_k$ and $G - P_{k+3}$ have $\Lambda $-factors, and so 
$G$ has a $\{\Lambda , P_k\}$- factor and a $\{\Lambda , P_{k + 3}\}$-factor.
\es
{\bf Proof} ~(uses {\bf \ref{frame,claw-free}}).
We prove $(a1)$ by induction on $v(G)$.
If $G$ is a cycle, then our claim is  obviously true.
Otherwise, consider $A = A(G)$ provided by Procedure 
${\cal E}$'. Then $e \not \in E(A)$.
Let $P$ be a maximum $\Lambda $-packing in $A(G)$.
Since $e(A) \ge 2$, clearly $v(P) = 3s$ for some $s \ge 1$.
Therefore $v(G - P) \equiv 0 \bmod 3$ and by  
{\bf \ref{frame,claw-free}}, $G - P$ is also a 2-connected claw-free graph. Since  $e \not \in E(A)$, clearly 
$e \not \in E(P)$. Obviously, $v(G - P) <  v(G)$.
By the induction hypothesis, $G - P$ has 
a $\Lambda $-factor $Q$ avoiding $e$. Then 
$P \cup Q$ is a $\Lambda $-factor of $G$ avoiding edge  $e$.
The proof of $(a2)$  is similar to that above. 
\ep
\\[1ex]
\indent
It turns out that an analogue of 
{\bf \ref{clfree,2con-avoid-e}} $(a2)$ when a 4-vertex path is replaced by a claw is also true provided a graph has a claw.

\bs
\label{G-Y}
Suppose that $G$ is a 2-connected  claw-free graph,
$v(G) \equiv 1 \bmod 3$, and $G$ is not a cycle. Then 
$G$ has at least two claws $Y$ such that $G - Y$ has a  
$\Lambda $-factor. 
\es

{\bf Proof} (uses {\bf \ref{Conclfr2endbplocks}} and 
{\bf \ref{frame,claw-free}}).
Let $H = F_c(G)$ and $A = A(G)$ (see  {\bf \ref{frame,claw-free}}), and so $A \subset H$. 
Suppose first that $A$ is a cycle. Then $A$ is a Hamiltonian cycle of $H$. Since $G$ is not a cycle, we have by {\bf \ref{frame,claw-free}} $(f2)$: $E(H) \setminus 
E(A) \ne \emptyset $ and every edge in $E(H) \setminus E(A)$ belongs to a unique triangle in $H$. Then $H$ has at least two claws and  $H - Y$ has a  $\Lambda $-factor
for every claw in $H$. 
Since $H$ is a spanning subgraph of $G$,  every  
$\Lambda $-factor of $H - Y$ is also a $\Lambda $-factor of $G - Y$. 

Now suppose that $A$ is a path. 
Let $x$ and $y$ be the end-vertices of $A$. 
Since $H$ is a spanning subgraph of $G$, it suffices  to prove the following
\\[0.5ex]
{\sc Claim.}
{\em 
For every vertex $v \in \{x,y\}$ there exist two claws $Y_v$ an $Z_v$ in $H$ such that 
either $G - Y_v$ or $G - Z_v$ has a $\Lambda $-factor.}
\\[0.5ex]
{\em Proof.}
By {\bf \ref{frame,claw-free}}, every end-vertex of $A$ has degree three and belongs to a unique triangle  of $H$.
By symmetry, we can assume that $v = y$.
Let $\Delta $ be the triangle containing $y$ and 
$V(\Delta ) = \{s,y,z\}$.
Let $Y$ and $Z$ be the claws in $H$ centered at $y$ and $z$, respectively. Then $Y$ contains the end-edge $yy'$ of $L$ and  $Z$ contains the edge $zz'$, where $z' \not \in V(\Delta )$.
Let $H' = H - A$. By {\bf \ref{frame,claw-free}} $(f3)$, $H'$ is a 2-connected claw-free spanning subgraph of $G - A$. 
\\[1ex]
${\bf (p1)}$
Suppose that $v(L)  \equiv 0 \bmod 3$.
Let $R_0 = (A- x) \cup Y $ and $H_0 = H - R_0$. 
Then  
$A' = R_0 - Y = A - \{x, y, y'\}$ is the subpath of $R_0$. Since $v(L)  \equiv 0 \bmod 3$, also 
$v(A')  \equiv 0 \bmod 3$. 
Then $A'$ has a unique $\Lambda $-factor $P$. 
Let $H'' = H - (A - x)$.
Obviously, $H_0 = H'' - \{s,z\}$.
Since $H'$ is 2-connected, $H''$ is also 2-connected. 
Since $s$ and $z$ have degree three in $H$, they both have degree two in $H''$. Therefore 
$H_0 = H'' - \{s,z\}$ has exactly two end-blocks. 
Since $H$ is claw-free, $H_0$ is also claw-free. 
By {\bf \ref{Conclfr2endbplocks}}, 
$H_0$ has a $\Lambda $-factor $Q$.
Then $P \cup Q$ is a $\Lambda $-factor of $H_0 \cup (A - Y) = H - Y$.
\\[1ex]
${\bf (p2)}$
Suppose that $v(A)  \equiv 1 \bmod 3$.
Let $R_0 = L \cup Z $ and $H_0 = H - R_0$.
Then  
$L' = R_0 - Z = L - y$ is the subpath of $R_0$. 
Since $v(A)  \equiv 1 \bmod 3$, clearly $v(A')  \equiv 0 \bmod 3$. 
Then $A'$ has a unique $\Lambda $-factor $P$. 
Since $s$ and $z$ have degree three in $H$, they both have degree two in $H'$. 
Now since $H'$ is 2-connected,  
$H_0 = H' - \{s,z,z'\}$ has exactly two end-blocks. 
Since $H$ is claw-free, $H_0$ is also claw-free. By {\bf \ref{Conclfr2endbplocks}}, 
$H_0$ has a $\Lambda $-factor $Q$.
Then $P \cup Q$ is a $\Lambda $-factor of $H_0 \cup (A - Z) = H - Z$.
\\[1ex]
${\bf (p3)}$
Suppose that $v(A)  \equiv 2 \bmod 3$.
Let $R_0 = L \cup \Delta $ and $H_0 = H - R_0$. 
Then 
$A' = R_0 - Y = A - \{y, y'\}$ is the subpath of $R_0$. 
Since $v(A)  \equiv 2 \bmod 3$, clearly $v(A')  \equiv 0 \bmod 3$. 
Then $A'$ has a unique $\Lambda $-factor $P$. 
Since $s$ and $z$ have degree three in $H$, they both have degree two in $H'$. 
Now since $H'$ is 2-connected,  
$H_0 = H' - \{s,z\}$ has exactly two end-blocks. 
Since $H$ is claw-free, $H_0$ is also claw-free. 
By {\bf \ref{Conclfr2endbplocks}}, 
$H_0$ has a $\Lambda $-factor $Q$.
Then $P \cup Q$ is a $\Lambda $-factor of $H_0 \cup (A - Y) = H - Y$.
\ep
\\[1ex]
\indent
By {\bf \ref{clfree,2con-avoid-e}} $(a2)$, every 
2-connected  claw-free graph with $v(G) \equiv 2 \bmod 3$
has an edge $xy$ such that $G - \{x,y\}$ has a  $\Lambda $-factor.
It turns that the following stronger result is  true.
\bs 
\label{clfree,2con,(x,y)-}
Suppose that $G$ is a 2-connected  claw-free graph and
$v(G) \equiv 2 \bmod 3$. Then  
for every vertex $x$ in $G$ there exist at least two edges $xb_1$ and $xb_2$ in $G$  such that each $G - \{x,b_i\}$ is connected and has a  $\Lambda $-factor.
\es

\bp (uses 
{\bf \ref{eb<3}},
{\bf \ref{reduction}}, and 
{\bf \ref{lambda(G)=}}).
Since $G$ is 2-connected, there exists 
an edge $xy$ in $G$ such that 
$G - \{x,y\}$ is connected.
Suppose that $G - \{x, y\}$  has no $\Lambda$ - factor.
Then by {\bf \ref{eb<3}},
$G -\{x,y\}$ has at least three end-blocks $B_i$, $i \in \{1, \ldots , k\}$, $k \ge 3$.
Let $b'_i$ be the boundary vertex of block $B_i$ in 
$G -\{x,y\}$. 
Let $V_i$ be the set of vertices in $\{x,y\}$ adjacent to a vertex in $B_i - b'_i$ and ${\cal B}_v$ be the set of the end-blocks in $G - \{x,y\}$ having an inner vertex adjacent to $v \in \{x,y\}$.
Since $G$ is 2-connected, each $|V_i| \ge 1$.
Since $G$ is claw-free, each $|{\cal B}_v| \le 2$.
Since $k \ge 3$, $|{\cal B}_z| = 2$ for some $z \in \{x,y\}$,
 say  ${\cal B}_z = \{B_1,B_2\}$.
Let $zb_i  \in E(G)$, where $b_i \in V(B_i - b'_i)$ for $i \in \{1,2\}$. 
Since $G$ is claw-free, $\{x, y, b_1, b_2\}$
does not induce a claw in $G$. Therefore $sb_j \in E(G)$ for $\{s,z\} = \{x,y\}$ and some $j \in \{1,2\}$, 
say, for $j = 2$. Then  ${\cal B}_s = \{B_i:  i \ge 2\}$. 
Since ${\cal B}_s = 2$, we have: $k = 3$ and
${\cal B}_s = \{B_2, B_ 3\}$.
Now we can assume that $z = x$ and $s = y$.
Obviously, $G - \{x,b_1\}$ is claw-free, connected, and has exactly two end-blocks. By  {\bf \ref{eb<3}}, 
$G - \{x,b_1\}$ has a $\Lambda $-factor.

We want to prove that $G - \{x,b_2\}$ also has a $\Lambda $-factor.
Let $C_i$ be the end-chain of $G - \{x,y\}$ containing 
$B_i$,  $i \in \{1,2,3\}$.
Graph $G - \{x,y\}$ is claw-free and has exactly three end-blocks.
Since $G - \{x,y\}$ has no $\Lambda$-factor, by 
{\bf \ref{lambda(G)=}} $(a3)$,
a graph 
$G_k$ obtained from $G$ by Reduction {\bf \ref{reduction}} has exactly three end-chains and each of them has one edge. Therefore each $v(C_i) \equiv 2 \bmod 3$.
Graph $G - \{x,b_2\}$ is  claw-free, connected,  and has a leaf $y$ and two or  three  end-chains. 
If $G - \{x,b_2\}$  has two  end-chains,
then by {\bf \ref{eb<3}}, $G - \{x,b_2\}$ has a $\Lambda$-factor.
So we assume that $G - \{x,b_2\}$  has three  end-chains
$C'_1$, $C'_1$, and $C'_1$, where 
$b'_1 \in V(C'_1)$ and $b'_3 \in V(C'_3)$.
Then $C'_1 = C_1$, $C'_2 = C_2 - b'_2$, and
$C'_3$ is obtained from $C_3$ by adding edge $yb_3$.
Since $v(C_2) \equiv 2 \bmod 3$, clearly $v(C'_2) \equiv 1 \bmod 3$.
Then a graph $G_k$ obtained from $G$ by Reduction
{\bf \ref{reduction}}  has two end-blocks.
Therefore by
{\bf \ref{lambda(G)=}} $(a3)$,
$G - \{x,b_2\}$ has a $\Lambda $-factor.
\ep
\\[2ex]
\indent
From {\bf \ref{clfree,2con,(x,y)-}} we have  for
3-connected claw-free graphs the following stronger result (with a simpler proof).
\bs 
\label{clfree,3con,(x,y)-} 
Suppose that $G$ is a 3-connected  claw-free 
graph and  $v(G) \equiv 2 \bmod 3$.
Then $G - \{x,y\}$ has a $\Lambda $-factor for every edge $xy$ in $G$.
\es

\bp  (uses {\bf \ref{eb<3}}). Let $G' = G - \{x,y\}$. Since $G$ is 3-connected, $G'$ is connected.
By {\bf \ref{eb<3}}, it suffices to
prove that $G'$ has at most two end-blocks.
Suppose, on the contrary, that $G'$ has at least three end-blocks. 
Let $B_i$, $i \in \{1,2,3\}$, be some three blocks of $G'$.
Since $G$ is 3-connected, for every block $B_i$ and every vertex $v \in \{x,y\}$ there is an edge $vb_i$, where $b_i$ is an inner vertex of $B_i$.  Then $\{v, b_1, b_2, b_3\}$  induces a claw in $G$, a contradiction.
\ep
\\[1ex]
\indent
As we have seen in the proof of  {\bf \ref{clfree,2con,(x,y)-}}, the claim of {\bf \ref{clfree,3con,(x,y)-}} is not true for claw-free graphs of connectivity two.
\bs
\label{clfree,3con,L-}
Suppose that $G$ is a 3-connected  claw-free graph and
$v(G) \equiv 0 \bmod 3$.
 Then for every edge $xy$ in $G$ there exist
at least two 3-vertex paths $L_1$ and $L_2$ in $G$ centered at $y$, containing $xy$, and such that each 
$G - L_i$ is connected and has a  $\Lambda $-factor.
\es

\bp (uses {\bf \ref{eb<3}}).
We need the following simple fact.
\\[1ex]
{\sc Claim}.
{\em Let $G$ be a 3-connected graph.  
 Then for every vertex $x$ and every edge $xy$ in $G$ there exist two 3-vertex paths $\Lambda _1$ and $\Lambda _2$ in $G$ centered at $y$, containing $xy$, and such that each 
 $G - \Lambda _i$ is connected.}
\\[1ex]
\indent
By the above  {\sc Claim}, $G$ has a 3-vertex path $L = xyz$ such that $G - L$ is connected. 
If every such 3-vertex path belongs to a $\Lambda $-factor of $G$, then we are done.
Therefore we assume that $G - L$ is connected but has no $\Lambda $-factor.
Obviously, $G - L$ is claw-free.
Therefore by {\bf \ref{eb<3}},
$G - L$ has at least three end-blocks $B_i$, 
$i \in \{1, \ldots , k\}$, $k \ge 3$.
Let $b'_i$ be the boundary vertex of $B_i$.
Let $V_i$ be the set of vertices in $L$ adjacent to inner 
vertices 
$G - L$ having an inner vertex 
adjacent to $v$ in $V(L)$.
Since $G$ is 3-connected, each $|V_i| \ge 2$.
 Since $G$ is claw-free, each $|{\cal B}_v| \le 2$.
It follows that $k = 3$, each $|V_i| = 2$, each $|{\cal B}_v| = 2$, as well as  all $V_i$'s are different and all ${\cal B}_v$'s are different. 
Let $s^1 = z$, $s^2 = x$,  $s^3 = y$, and 
$S = \{s^1, s^2, s^3\}$. 
We can assume that $V_i = S - s^i$, $i \in \{1,2,3\}$.
Then for every vertex $s^j \in V_i$ there is a vertex  $b_i^j$ in  $B_i - b'_i$ adjacent to $s^j$, where $\{b_i^j: s^j \in V_i\}$ has exactly one vertex if and only if $B_i - b'_i$ has exactly one vertex.
Let $L_i = s^2s^3b_i$, where $b_i = b_i^3$.
By {\bf \ref{eb<3}},
it suffices to show that each $G - L_i$ is connected and 
has  at most two end-blocks. 

Let $i = 1$. If $B_1 - b_1$ is 2-connected, then 
$B_1 - b_1$ and $G - L_1 - (B_1 - b'_1)$ are the two 
end-blocks of $G - L_1$ and we are done. 
If $B_1 - b_1$ is empty, then $G - L_1$ is 2-connected.
So we assume that
$B_1 - b_1$ is not empty and not 2-connected. Then 
$B_1 - b_1$ is connected and has exactly two end-blocks, say $C_1$ and $C_2$. Let $c'_i$ be the boundary vertex of $C_i$ in $B_1 - b_1$. Since $G$ is 3-connected, each 
$C_i -c'_i$ has a vertex  adjacent to $\{s^2,s^3\}$.
We can assume that a vertex $c_1$ in $C_1 - c'_1$ is adjacent to $s^2$.
If there exists a vertex $c_2$  in $C_2 - c'_2$ adjacent to 
$s^2$,
then 
$\{s^2, b_3^2, c_1, c_2\}$ 
induces a claw in $G$, a contradiction. So we assume that  no vertex in $C_2 - c'_2$ is adjacent to $s^2$.
Then there is a vertex $c_2$ in $C_2 - c'_2$ adjacent to 
$s^3$.
Then $\{s^2, s^3, b_2^3, c_2\}$
 induces a claw in $G$, a contradiction.

Finally, let $i = 2$.
If $B_2 - b_2$ is 2-connected, then 
$B_1$ and $G - L_2 - (B_1 - b'_1)$ are the two 
end-blocks of $G - L_2$ and we are done. If $B_1 - b_1$ is empty, then $G - L_2$ has two end-blocks, namely, $B_1$ and the subgraph of $G$ induced by $B_3 \cup s^1$.
So we assume that
$B_2 - b_2$ is not empty and not 2-connected. Then 
$B_2 - b_2$ is connected and has exactly two end-blocks, say $D_1$ and $D_2$. Let $d'_i$ be the boundary vertex of $D_i$ in $B_2 - b_2$. Since $G$ is 3-connected, each 
$D_i -d'_i$ has a  vertex  adjacent to $\{s^1,s^3\}$.
We can assume that a vertex $d_1$ in $D_1 - d'_1$ is adjacent to $s^3$.
If there exists a vertex $d_2$  in $D_2 - d'_2$ adjacent to 
$s^3$,  
then $\{s^3,d_1,d_2, b_1^3\}$ induces a claw in $G$, a contradiction.
So suppose that no vertex in $D_2 - d'_2$ is adjacent to 
$s^3$.
Then there is a vertex $d_2$ in $D_2 - d'_2$ adjacent to 
$s^1$.
Then $\{s^1, s^3, b_3^1, d_2\}$  
induces a claw in $G$, a contradiction.
\ep
\\[2ex]
\indent
From the proof of {\bf \ref{clfree,3con,L-}} we  have, in particular:
\bs 
\label{clfree,3con,deg3cntrL-}
Suppose that $G$ is a 3-connected  claw-free graph and 
$v(G) \equiv 0 \bmod 3$.
If $L$ is a 3-vertex path and the center vertex of $L$ has degree 3 in $G$, then  $G - L$ is connected and has 
a  $\Lambda $-factor in $G$.
\es

From {\bf \ref{clfree,3con,deg3cntrL-}} and the proof of 
{\bf \ref{clfree,3con,L-}} we  have:
\bs 
\label{clfree,4con,deg3cntrL-}
Suppose that $G$ is a cubic 3-connected claw-free graph or 4-connected  claw-free graph
 with $v(G) \equiv 0 \bmod 3$.
Then $G - L$ is connected and has a  $\Lambda $-factor for every 3-vertex path $L$ in $G$.
\es

The claim of {\bf \ref{clfree,4con,deg3cntrL-}} may not be true for a  claw-free graph of connectivity 3 if they are not cubic.
Recall that a {\em net} is a graph obtained from a claw  by replacing its vertex of degree 3 by a triangle.
Let $N$ be a net with the three leaves $v_1$, $v_2$, and $v_3$,
$T$  a triangle with $V(T) = \{t_1, t_2, t_3\}$, and let $N$ and $T$ be disjoint. 
Let $H = N \cup T \cup \{v_it_j: i, j \in \{1,2,3\}, i \ne j\}$.
Then $H$ is a 3-connected claw-free graph, $v(H) = 9$, each $d(t_i, H) = 4$, $d(x, H) = 3$ for every $x \in V(H - T)$, and  
$ H - T = N$ has no $\Lambda $-factor.  If $L$ is a 3-vertex path in $T$, then 
$H - L = H - T$, and so $H - L$ has no $\Lambda $-factor.
There are infinitely many pairs $(G, L)$ such that
$G$ is a 3-connected, claw-free, and non-cubic graph, 
$v(G) \equiv 0 \bmod 3$, $L$ is a 3-vertex path in $G$, and 
$G - L$ has no 
$\Lambda $-factor. 

By {\bf \ref{A}}, such a pair can be obtained from the above pair $(H,L)$ by replacing $N$ by any graph $A$ with
exactly three leaves satisfying the assumptions of
 {\bf \ref{A}}.
\\[2ex]
\indent
From {\bf \ref{clfree,3con,L-}} we have, in particular:
\bs 
\label{clfree,3con,einL}
Suppose that $G$ is a 3-connected  claw-free graph
and $v(G) \equiv 0 \bmod 3$.
Then for every edge $e$ of $G$
there exists a $\Lambda $-factor in $G$ containing $e$.
\es

The following examples show that assumption
``{\em $G$ is a 3-connected graph}'' in 
{\bf \ref{clfree,3con,einL}} is essential.
Let $R$ be the graph obtained from two disjoint cycles $A$ and $B$ by adding a new vertex $z$,  and the set of four new edges $\{a_iz, b_iz: i \in \{1,2\}\}$, 
where $a = a_1a_2 \in E(A)$ and $b = b_1b_2 \in E(B)$. 
It is easy to see that $R$ is a claw-free graph of connectivity one. 
Furthermore, if $v(A) \equiv 1 \bmod 3$ and 
$v(B) \equiv 1 \bmod 3$, then $v(R) \equiv 0 \bmod 3$ and 
$R$ has no 
$\Lambda $-factor containing edge $e \in \{a, b\}$.
Similarly,
let $Q$ be the graph obtained from two disjoint cycles $A$ and $B$ by adding two new vertices $z_1$ and $z_2$, a new edge $e = z_1z_2$,  and the set of eight new edges $\{a_iz_j, b_iz_j: i,j \in \{1,2\}\}$, where $a_1a_2 \in E(A)$ and $b_1b_2 \in E(B)$. 
It is easy to see that $Q$ is a claw-free graph of connectivity two. Furthermore, if $v(A) \equiv 2 \bmod 3$ and 
$v(B) \equiv 2 \bmod 3$, then $v(Q) \equiv 0 \bmod 3$ and $Q$ has no $\Lambda $-factor containing edge $e$.
\\[1ex]
\indent
Let $F$ be a graph, $x \in V(F)$, and $X = \{x_1,x_2, x_3\}$ be the set of  vertices in $F$ adjacent to $x$.
Let $T$ be a triangle, $V(T) = \{t_1, t_2, t_3\}$, and $V(F) \cap V(T) = \emptyset $. 
Let $G = (F - x) \cup T \cup \{x_it_i: i \in \{1,2,3\}\}$.
We say that $G$ {\em is obtained from $F$ by replacing a vertex $x$ by a triangle}. 

Given a cubic graph $F$ with possible parallel edges, 
let $F^\Delta $ denote the graph obtained from  $F$ by replacing each 
vertex  of $F$ by a triangle. Clearly, $F^\Delta $ is cubic and claw-free, every vertex  belongs to exactly one triangle, every edge belongs to at most one triangle in $F^\Delta $, and $v(F^\Delta ) \equiv 0 \bmod 3$. 
Obviously, $F^\Delta $ is $k$-connected if and only if 
$F$ is $k$-connected, $k \in \{1,2,3\}$.
We call $F^\Delta $ a {\em $\Delta $-graph}.
\bs
\label{2con,cubic,tr}
Let $G$ be a 2-connected  $\Delta $-graph. 
 Let $L$ be a 3-vertex path in $G$.
Then
\\[0.5ex]
$(a)$ $G - L$ has a $\Lambda$-factor.

Moreover,
\\[0.5ex]
$(a1)$ if $L$ induces  a triangle in $G$, then 
$G$ has a $\Lambda$-factor $R$ containing $L$ and such that each component of $R$ induces a triangle
\\[0.5ex]
$(a2)$ if $L$ does not induce a triangle in $G$, then 
$G$ has a $\Lambda $-factor $R$ containing $L$ and 
such that no component of $R$ induces a triangle, and
\\[0.5ex]
$(a3)$ if $L$ does not induce a triangle in $G$, 
then $G$ has  a $\Lambda $-factor 
containing $L$ and a component that induces a triangle.
\es

\bp Since $G$ is a 2-connected  $\Delta $-graph, $G$ can be obtained from  a 2-connected cubic graph $G'$ (with possible parallel edges) by replacing each vertex of $G'$ by a triangle.
Obviously, there is a natural bijection $\alpha : E(G') \to E'$. Let $E'$ be the set of edges in $G$ that
belong to no triangle. Let $L = xzz_1$.
Since each vertex of $G$ belongs to exactly one triangle,
we can assume that  $xz$ belongs to a triangle $T = xzs$.
\\[1ex]
${\bf (p1)}$ Suppose that $L$ induces a triangle in $G$, 
and so $s = z_1$.
Obviously the union of all triangles in $G$ contains
a $\Lambda $-factor, say $P$, of $G$ and 
$L \subset P$.
Therefore claim $(a1)$ is true.
\\[1ex]
${\bf (p2)}$ Now suppose that $L$ does not induce a triangle in $G$, and so $s \ne z_1$. 
Let $\bar{s} = ss_1$ and $\bar{z} = zz_1$ be the edges of $G$ not belonging to $T$, and therefore belonging to no triangles in $G$. Hence $\bar{s} = \alpha (\bar{s}')$ and 
$\bar{z} = \alpha (\bar{z}')$, where
$\bar{s}' = s's'_1$ and $\bar{z}' = z'z'_1$ are edges in $G'$, and $s'  = z'$.
Since every vertex in $G$ belongs to exactly one triangle, clearly $s_1 \ne z_1$.
\\[1ex]
${\bf (p2.1)}$
We prove $(a2)$. 
By using Tutte's criterion for a graph to have a perfect matching (see \cite{LP}), it is easy to prove the following 
\\[1ex]
{\sc Claim}.  
{\em If $A$ is a cubic 2-connected
graph, then for every 3-vertex path $J$ of $A$ there exists 
a 2-factor of $A$ containing $J$.}
 \\[1ex]
 \indent
By  the above {\sc Claim}, $G'$ has a 2-factor $F'$ containing 3-vertex path $S' = s'_1s'z'_1$.
Let $C'$ be the (cycle) component of $F'$ containing $S'$.
If $Q'$ is a (cycle) component of $F'$, then let $Q$ 
be the subgraph of $G$, induced by the edge subset 
$\{\alpha (e): e \in E(Q')\} \cup \{E(\Delta _v): v \in V(Q')\}$.
Obviously $v(Q) \equiv 0 \bmod 3$ and $Q$ has a (unique) Hamiltonian cycle $H(Q)$. 
Also the union $F$ of all $Q$'s is a spanning subgraph of 
$G$ and each $Q$ is a component of $F$.  
Moreover, if $C$ is the component in $F$, corresponding to $C'$, then $L \subset H(C)$.
Therefore each $H(Q)$ has a  $\Lambda $-factor $P(Q)$, such that no component of $P(Q)$ induces a triangle, and 
$H(C)$ has a (unique) $\Lambda $-factor $P(C)$, such that $L \subset P(C)$ and no component of $P(C)$ induces 
a triangle.
The union of all these $\Lambda $-factors
is a $\Lambda $-factor $P$ of $G$ containing $L$ and 
such that no component of $P$ induces a triangle.
Therefore $(a2)$ holds.
\\[1ex]
${\bf (p2.2)}$
Finally, we prove $(a3)$.
Since $G'$ is 2-connected and cubic, there is a cycle $C'$ in $G'$ such that $V(C') \ne V(G')$ and $C'$ 
contains $S' = s'_1s'z'_1$.
Let, as above,  $C$ be the subgraph of $G$, induced by 
the edge subset
$\{\alpha (e): e \in E(C')\} \cup \{E(\Delta _v): v \in V(C')\}$.
Obviously, $v(C) \equiv 0 \bmod 3$, $C$ has a (unique) Hamiltonian cycle $H$, and $L \subset H$.
Therefore $H$ has a (unique) $\Lambda $-factor $P(C)$ containing $L$.
Since $V(C') \ne V(G')$, we have $V(G' - C') \ne \emptyset $. Therefore $G - C$ has a triangle. Moreover, every vertex $v$ in $G - C$ belongs to a unique triangle $\Delta _v$, and therefore as in ${\bf (p1)}$, $G - C$ has a $\Lambda $-factor $Q$ whose every component induces a triangle in $G - C$. Then $P(C) \cup Q$ is a required a $\Lambda $-factor  in $G$.
\ep 
\\[1ex]
\indent
Obviously, {\bf \ref{2con,cubic,tr}} $(a)$ also follows from 
{\bf \ref{clfree,3con,deg3cntrL-}}.

Theorem {\bf \ref{2con,cubic,tr}} is not true for a cubic 2-connected claw-free graph $F$  with an edge $xy$ belonging to two triangles $T_i$ with $V(T_i) = \{x,y, z_i\}$, $i \in \{1,2\}$, because $L = z_1xz_2$ is a 3-vertex path in $F$ and $y$ is an isolated vertex in $F - L$.
\\[1ex]
\indent
Now we can give polynomial-time characterization of pairs $(G,E)$ such that $G$ is a 2-connected  $\Delta $-graph, $E \subset E(G)$, $|E| = 3$, and $G - E$ has no $\Lambda $-factor.
Recall that if $E \subseteq E(G)$, then $\dot{E}$ denotes the subgraph of $G$ induced by $E$.

\bs 
\label{clfree,2con-avoid(a,b,c)}
Suppose that $G$ is a 2-connected  $\Delta $-graph.
Let  $E \subset E(G)$ and  $|E| = 3$. Then 
 the following are equivalent:
\\[.5ex]
$(g)$
$G - E$ has no $\Lambda $-factor and
\\[.5ex]
$(e)$ $\dot{E}$ satisfies one of the following conditions:
\\[.5ex]
\indent
$(e1)$ $\dot{E}$ is  a claw,
\\[.5ex]
\indent
$(e2)$ $\dot{E}$ is a triangle,
\\[.5ex]
\indent
$(e3)$
$\dot{E}$ has exactly two components, 
the 2-edge component $\dot{E}_2$  belongs to a triangle in $G$,
the 1-edge component $\dot{E}_1$ belongs to no triangle in $G$, and $G - E$ is not connected, and 
\\[.5ex]
\indent
$(e4)$
$\dot{E}$ has exactly two components,  
the 2-edge component  $\dot{E}_2$ belongs to a triangle $T$ and the 1-edge component $\dot{E}_1$  belongs to a triangle $D$ in $G$,  $\dot{E}_1$ and $\dot{E}_2$ belong to different component of 
$G - \{ d,t \}$, where
$d$ and $t$ are the edges in $G - E(D) - E(T)$ incident to the single vertex of $D - \dot{E}_1$ and to the isolated vertex of $T - E_2$, respectively.
 \es

\bp (uses {\bf \ref{eb<3}},
{\bf \ref{A}}, and
{\bf \ref{2con,cubic,tr}}(a)). 
Let $E = \{a,b,c\}$, where $c = c_1c_2$. 
\\[1ex]
${\bf (p1)}$ We prove $(e) \Rightarrow (g)$.

Suppose that $\dot{E}$ satisfies $(e1)$, i.e. $\dot{E}$ is a claw.
Then $G - E$ has an isolated vertex and therefore has  no $\Lambda $-factor. 

Suppose that $\dot{E}$ satisfies $(e2)$, i.e. $\dot{E}$ is a triangle. 
Then by {\bf \ref{A}}, $G - E$ has no 
$\Lambda $-factor.

Now we  assume (as in $(e3)$ and $(e4)$) that
$\dot{E}$ has exactly two components $\dot{E}_2$ and 
$\dot{E}_1$,
the 2-edge component $\dot{E}_2$  belongs to a triangle 
$T$ in $G$, $E_2 = \{a,b\}$, and $E_1 = \{c\}$,
$t = t_1t_2$ is an edge in $G - E$, where $t_1$ is an isolated vertex in $T - E_2$, and so $t_1$ is a (unique) leaf
in $G - E$.
Let $u$ be the edge in $T$ distinct from $a$ and $b$.

Suppose that $\dot{E}$ satisfies $(e3)$,
and so $\dot{E}_1$ belongs to no triangle in $G$ and
$G - E$ is not connected. Obviously, $G - E$ has exactly two components. Let $S$  the component in $G - E$ containing edge $t$. Then edge $u$ is not in $S$.
Therefore every vertex in $S$ distinct from the leaf $t_1$  belongs to exactly one triangle. 
Hence $v(S) \equiv 1\bmod 3$ implying that $G - E$ has no $\Lambda $-factor.

Finally, suppose that $\dot{E}$ satisfies $(e4)$, and so
edge $c = c_1c_2$ belongs to a triangle $D$ in $G$,
$G - \{d,t\}$ is not connected, and $E_1$, $E_2$ belong to different components of $G - \{d,t\}$, where $d$ is the edge in $G - E(D)$ incident to the single vertex in $D - \{c_1,c_2\}$.
Suppose, on the contrary, that $G - E$ has a 
$\Lambda $-factor $P$.  
Since $G$ is 2-connected and claw-free, $G - \{d,t\}$ has exactly two components. 
Therefore 
$G - \{a,b\} - (E(D) - c)$ has also two components.
Let $C'$ be the component of $G - \{a,b\} - (E(D) - c)$ containing $c$.
Since $G$ is a $\Delta $-graph, $G = F^{\Delta }$ for some cubic 2-connected graph $F$ (with possible parallel edges). Let $d'$ and $t'$ be the edges in 
$F$ corresponding to edges  $d$ and $t$ of $G$, respectively. 
Since $F$ is 2-connected,  $F - \{d',t'\}$ has at most two components.  Since $G - \{d,t\}$ is not connected, 
$F - \{d',t'\}$ has exactly two components.
 It follows that
$H = G - \{a,b\} - E(D)) = G - E - E(D)$ has also exactly two components.
Therefore $C = C' - c$ is connected, and so $C$ is a component of $H$ containing the end-vertices $c_1$ and $c_2$ of edge $c$.

Let $C_u$ and $C_t$ be the components of $H$ containing $u$ and $t$, respectively. Then $C_u \ne C_t$.
 Now $C_2 = C_u \cup T$ is the component in $G - \{d,t\}$ containing $E_2$. By $(e4)$, $c \not \in E(C_2)$.
 Therefore $c$ is an edge of   $C_1 = C_t \cup D$. 
 Thus $C_t = C$.
 Clearly, $C$ has exactly three leaves $c_1$, $c_2$, and $t_1$ (the leaf incident to $t$) and every other vertex of 
 $C$ belongs to a unique triangle in $C$, and so $v(C) \equiv 0 \bmod 3$. 
 By {\bf \ref{A}}, $C$ has no 
$\Lambda $-factor. Therefore $P$ has a 3-vertex path $L$ which contains at least one edge in $D - c$. 
Since $t$ is a dangling edge in $G - E$, clearly $P$ also has a 3-vertex path $L_t$ containing $t$ and 
$L_t \subset C$. Therefore $P$ has to contain a 
$\Lambda $-factor of $C - L - L_t$. 
However, $v(C - L - L_t) \not \equiv 0 \bmod 3$, a contradiction.
\\[1ex]
${\bf (p2)}$ Finally, we prove $(g) \Rightarrow (e)$.
Namely, we assume that $\dot {E}$ does not satisfy $(e)$ and we want to show that in this case $G - E$ has a 
$\Lambda $-factor. 

Let $X, Y \subset E(G)$ be such that $X$ meets no triangle in $G$, each edge in $Y$ belongs to a triangle in $G$,  and  no triangle in $G$ has more than one edge from $Y$, and so  $X \cap Y = \emptyset $.
We will use the following simple observation.
\\[0.5ex]
{\sc Claim.}
{\em 
$G - X - Y$ has a  $\Lambda $-factor $P$ such that  every component of $P$ induces a triangle in $G$ and if an edge $y$ from $Y$ is in a triangle $T$, then $T - y$ is a component of $P$.}
\\[0.5ex]
\indent
By the above {\sc Claim}, we can assume that the two edges of $E_2$ belong to the same triangle $T$. 

Suppose that $\dot {E}$ is connected.
Since $\dot {E}$ does not satisfy $(e)$, $\dot {E}$ is not a claw and not a triangle. Then $\dot {E}$ is a 3-edge path and  $u, t \not \in E$.
Let $V$ be a 3-vertex path in $G$ containing $u$ and avoiding $E$. 
Then $G - V$ has no edges from $E$,  and so 
$G - V = G - E -V$.
By {\bf \ref{2con,cubic,tr}}(a), $G - V$ has a $\Lambda $-factor.

Finally, suppose that $\dot {E}$ is not connected, and so
$\dot {E}$ has exactly two components $\dot {E}_1$ and 
$\dot {E}_2$.
As in $\bf (p1)$, let $E_2 = \{a,b\}$ and $E_1 = \{c\}$, and
let $u$ be the edge of $T$ distinct from $a$ and $b$.
\\[1ex]
${\bf (p2.1)}$ Suppose that $c$ belongs to no triangle in $G$.
Since $\dot {E}$ does not satisfy $(e)$ (namely, $(e3)$), 
$G - E$ is connected. 
Clearly, $G - E$ is claw-free. 
Also $G - E$ has exactly two end-blocks and the  block of one edge $t$ is one of them.
By {\bf \ref{eb<3}}, $G - E$ has a $\Lambda $-factor.
\\[1ex]
${\bf (p2.2)}$ 
Finally, suppose that $c$ belongs to a triangle $D$ in $G$. Then $D \ne T$. 
Since $\dot {E}$ does not satisfy $(e)$ (namely, $(e4)$), 
$\dot {E}_1$ and $\dot {E}_2$ belong to the same component of $G - \{d,t\}$.
Let $V(D) = \{c_1, c_2, d_1\}$ and 
as above $c = c_1c_2$. Let $d = d_1d_2$ and $t = t_1t_2$, where $t_1 \in V(T)$, and so $t_1$ is an isolated vertex in $T - \{a,b\}$.

Let $G' = G - \{c,d\}$. Then $G'$ and $G' - E_2$ are claw-free and $v(B) \equiv 0 \bmod 3$ for every block $B$ in $G'$. Obviously, $G' - E_2 = G - E - d$.
Since $G$ is 2-connected, $G - c$ is also 2-connected. Therefore $eb(G') \le 2$ and if $eb(G') = 2$, then 
the end-vertices $d_1$ and $d_2$ of edge $d$ belong to different end-blocks $B_1$ and $B_2$ of $G'$, respectively.
If $eb(G') = 1$, then $G'$ is 2-connected. Therefore 
$G' - E_2$ has at most two end-blocks. By 
{\bf \ref{eb<3}},
$G' - E_2  = G - E - d$ has a $\Lambda $-factor  $P$, which is also a  $\Lambda $-factor of $G - E$.
So we assume that $eb(G') = 2$.
Since $E_1$ and $E_2$ belong to the same component of 
$G - \{d,t\}$,  clearly $d_1$, and $E_2$ belong to the same component of $G' - t$.

Suppose that $t$ is a cut-edge of  $G'$. 
Then $G' - t$ has exactly two components $C_1$ and
 $C_2$ containing $\{d_1,t_1\}$ and $d_2$, respectively, and each $v(C_i) \equiv 0 \bmod 3$. 
Let $L$ be a 3-vertex path $c_1d_1d_2$ in $G - c$. 
Then $C'_1 = C_1 - \{c_1, d_1, t_1\}$ and
$C'_2 = (C_2 - d_2) \cup \{t,t_1\}$ are the two components of $G - E - L$ containing $d_1$ and $d_2$,  respectively. Also each, $v(C'_i) \equiv 0 \bmod 3$ and $C'_i$ is claw-free and has exactly two end-blocks. 
By {\bf \ref{2conclfr}} and  {\bf \ref{Conclfr2endbplocks}}, $C'_i$ has a $\Lambda $-factor  $P_i$, 
$i \in \{1,2\}$. Then $L \cup P_1 \cup P_2$ is 
a $\Lambda $-factor  of $G - E$.

Finally, suppose that  $t$ is not a cut-edge of  $G'$.
Then $\{a,b,t\}$ belongs to a 2-connected block $R$ of $G'$. Then $R - E_2$ has at most two end-blocks. 
By {\bf \ref{2conclfr}} and  {\bf \ref{Conclfr2endbplocks}}, 
$R - E_2$ has a $\Lambda $-factor $Q$.
Let ${\cal B}$ denote the set of all 2-connected blocks of $G'$ distinct from $R$. Since $v(B) \equiv 0\bmod 3$ for every $B \in {\cal B}$,
each $B$ in ${\cal B}$ has a $\Lambda $-factor $P(B)$.
Then $\{P(B): B \in {\cal B}\} \cup Q$ is a $\Lambda $-factor of $G'$, which is also a $\Lambda $-factor of $G - E$.
\ep
\\[1ex]
\indent
From {\bf \ref{clfree,2con-avoid(a,b,c)}} we have, in particular:
\bs 
\label{clfree,2con-avoid(a,b)}
Suppose that $G$ is a $\Delta $-graph.
Let  $E \subset E(G)$ and  $|E| = 2$. Then
$G - E$ has a $\Lambda $-factor.
\es

From {\bf \ref{clfree,2con-avoid(a,b,c)}} we also have:
\bs 
\label{clfree,3con-avoid(a,b,c)}
Suppose that  $G$ is a 3-connected claw-free graph.
Let  $E \subset E(G)$ and  $|E| = 3$. Then
$G - E$ has a $\Lambda $-factor if and only if 
$\dot {E}$ is not a claw and not a triangle.
\es

It turns out that condition ``{\em $G$ is claw-free}'' in  {\bf \ref{clfree,2con-avoid(a,b,c)}} and in {\bf \ref{clfree,3con-avoid(a,b,c)}} is essential.
Namely, we have a construction showing that for every 3-edge graph $Y$ with no isolated vertices there are infinitely many pairs $(G, E)$ such that
$G$ is a cubic 3-connected graph, $v(G) \equiv 0 \bmod 3$, 
$E \subset E(G)$,  the subgraph $\dot{E}$ induced by $E$ in $G$ is isomorphic to $Y$ (and so $|E| = 3$), and 
$G - E$ has no $\Lambda $-factor.
\bs 
\label{clfree,2con,x-}
Suppose that $G$ is a 2-connected  claw-free graph and 
$v(G) \equiv 1 \bmod 3$.  
Then $G - x$ has a $\Lambda $-factor for every vertex $x$ in $G$.
\es

\bp (uses {\bf \ref{eb<3}}).
Let $x \in V(G)$. 
Since $v(G) \equiv 1 \bmod 3$, clearly $v(G - x) \equiv 0 \bmod 3$.
Since $G$ is 2-connected, $G - x$ is connected.
Since $G$ is  claw-free, $G - x$ is claw-free and has 
at most two end-blocks.
By {\bf \ref{eb<3}},
$G - x$ has a $\Lambda $-factor.
\ep
\\[1ex]
\indent
If on step $(s1)$ of Procedure ${\cal E}$ we find a longest cycle in $G$ containing a given vertex $x$, then this modification of Procedure ${\cal E}$ can also be used to prove {\bf \ref{clfree,2con,x-}}.
\\[1ex]
\indent
Moreover, the following strengthening of 
{\bf \ref{clfree,2con,x-}} holds for 3-connected  claw-free  graphs.
\bs 
\label{clfree,3con,(x,e)-} 
Suppose that $G$ is a 3-connected  claw-free 
graph and  $v(G) \equiv 1 \bmod 3$. Then $G - \{x,e\}$ has a $\Lambda $-factor for every vertex $x$ and every edge $e$ in $G$.
\es

\bp (uses {\bf \ref{2conclfr}} and {\bf \ref{clfree,2con-avoid-e}}).
Since $G$ is 3-connected, $G - x$ is a 2-connected  claw-free graph. Since $v(G) \equiv 1 \bmod 3$, we have
$v(G - x) \equiv 0 \bmod 3$. 
By {\bf \ref{2conclfr}}, $G - x$ has a $\Lambda $-factor $P$.
If $e \not \in E(G - x)$, then $P$ is a $\Lambda $-factor  of 
$G - \{x,e\}$. If $e  \in E(G - x)$, then by 
{\bf \ref{clfree,2con-avoid-e}}, $G - \{x,e\}$ has 
a $\Lambda $-factor.
\ep

\section{Packings and domination in graphs}

\label{domination}

\indent

Recall that $X \subseteq V(G)$ is called a domination set in graph $G$, if every vertex in $V(G) \setminus X$ is  adjacent to a vertex in $X$ and that the  domination number
$\gamma (G)$ is  the size of a minimum  
domination set in $G$.
We call a subgraph $P$ in $G$ a {\em star-packing} if every component of $P$ is isomorphic to $K_{1,s}$ for some integer $s \ge 0$. 

Obviously, $X$ is a domination set in $G$ if and only if there exists a star-factor $P = P(X)$ such that $Cmp(P) = \{P_x: x \in X\}$, where  
$x \in V(P_x)$ and $x$ is a (unique) vertex of degree at least two if $v(P_x) \ge 3$, and so $|Cmp(P)| = |X|$. 
Thus $X$ is a minimum domination set in $G$ if and only if 
$P(X)$ is a star-factor in $G$ having the minimum number of components and $\gamma (G) = cmp (P(X))$.

It is easy to show that every connected graph $G$ with no isolated  vertices  has a star-factor with no isolated vertices, and so $\gamma (G) \le v(G)/2$.

Clearly, every $\Lambda $-packing $P$ is a star-packing in $G$  and 
$P$ can be extended to a star-factor $P'$ in $G$. Then 
$\gamma (G) \le cmp(P')$.
For that reason,  results on the maximum $\Lambda $-packings in graphs may be useful in the study of some  graph domination problems.
Here is an example of such correlation.

In \cite{R} B. Reed conjectured that
if $G$ is a connected cubic graph, then 
$\gamma (G) \le \lceil v(G)/3 \rceil $.
It turns out that Reed's conjecture is not true for connected and even for 2-connected cubic graphs \cite{Kcntrex,KS}. 
Obviously,
\\[0.7ex]
$(d1)$
if a graph $G$ has 
a $\Lambda $-factor (and so $v(G) \equiv 0 \bmod 3$), then 
$\gamma (G) \le v(G)/3$,
\\[0.7ex]
$(d2)$
if $v(G) \equiv 1 \bmod 3$ and $G - x$ has a $\Lambda $-factor
for some vertex $x$ of $G$,
then $\gamma (G) \le \lceil v(G)/3 \rceil $, and
\\[0.7ex]
$(d3)$
if $v(G) \equiv 2 \bmod 3$ and $G - \{x,y\}$ has 
a $\Lambda $-factor
for some edge $xy$ of $G$,
then again $\gamma (G) \le \lceil v(G)/3 \rceil $.
\\[0.7ex]
\indent
Now if claim $(P)$ in Problem {\bf \ref{Pr3con}}  is true, then from {\bf \ref{cubic3-conZTF}} and {\bf \ref{cubic3-conT}}
it follows, in particular, that $(d1)$, $(d2)$, and $(d3)$ above are true, and so Reed's conjecture is  true for 3-connected cubic graphs.
\\[1ex]
\indent
The following packing result is also related with Reed's domination conjecture.
\bs {\em \cite{KcubHam}}
\label{Ham}
If $G$ is a cubic Hamiltonian graph with 
$v(G) \equiv 1\bmod 3$, then $G$ has a claw $Y$ 
such that
$G - Y$ has a $\Lambda $-factor, and so $G$ has 
$\{\Lambda , Y\}$-factor. 
\es

It follows that if $G$ is a cubic Hamiltonian graph with 
$v(G) \equiv 1\bmod 3$, then 
$\gamma (G) \le \lfloor v(G)/3 \rfloor $ which is stronger than Reed's conjecture suggests. 
The following natural question arises:
\bs {\em {\bf Problem}}
\label{problem,v(G)=1mod3}
Is it true that $\gamma (G) \le \lfloor v(G)/3 \rfloor $
for every  cubic 3-connected graph $G$ with
$v(G) \equiv 1\bmod 3$ ?
\es

In \cite{Kcntrex} we gave a construction providing infinitely many cubic cyclically 4-connected graphs 
$G$ with $v(G) \in \{0,2\} \bmod 3$ for which 
$\gamma (G) = \lceil v(G)/3 \rceil $,  and so Reed's suggested bound is tight even in the class of cyclically 
4-connected graphs.
From this construction it also follows that
the claim similar to {\bf \ref{Ham}} for graphs $G$ with 
$v(G) \not \equiv 1 \bmod 3$ is not true, namely, the graphs  provided by the construction have no 
$\{\Lambda , Y\}$-factor.

No bipartite counterexamples to Reed's conjecture have been found. We can show that if claim $(P)$ in Problem {\bf \ref{Pr3con}}  is true for bipartite graphs, then Reed's conjecture is  also true for bipartite 3-connected cubic graphs. 
\\[1ex]
\indent
Let $\gamma _i(G)$ denote the size of a minimum independent domination set, and so $\gamma _i(G) \ge \gamma (G)$.
It is easy to see that if $G$ is a claw-free graph, then
$\gamma _i(G) = \gamma (G)$.
\\[1ex]
\indent
From {\bf \ref{clfree,2con-avoid-e}} and {\bf \ref{G-Y}} we have the following upper bounds on the domination number of claw-free graphs:
\bs
\label{gamma}
Let $G$ be a 2-connected  claw-free graph.
Then  $\gamma (G) \le \lceil v(G)/3 \rceil $ and if, in addition, $G$ is not a cycle and $v(G) \equiv 1 \bmod 3$, 
then $\gamma (G) = \gamma _i(G) \le \lfloor v(G)/3 \rfloor $.
\es

\section{Further related results and questions}
\label{related-results}



\indent

Given a family ${\cal F}$ of non-isomorphic graphs,
an {\em edge disjoint ${\cal F}$-packing}
${\cal Q}$ of $G$ is a set $\{ Q_1, \ldots Q_k \}$ such that each 
$Q_i \subseteq E(G)$, every two members of ${\cal Q}$ are disjoint, and the subgraph $\dot{Q}_i$ induced by $Q_i$ in $G$ is isomorphic to a member of ${\cal F}$. 
Let  $E({\cal Q}) = \cup \{ E(Q_i): i \in \{1, \ldots k \} \}$
and $k({\cal Q}) = k$.  
 An edge disjoint ${\cal F}$-packing ${\cal Q}$ 
in $G$ is called 
an {\em edge disjoint ${\cal F}$-factor} of $G$ if $E(P) = E(G)$. 
The {\em edge disjoint ${\cal F}$-packing problem}
is the problem of 
finding in $G$ an edge ${\cal F}$-packing ${\cal Q}$ having the maximum number of edges $|E({\cal Q})|$.
If ${\cal F}$ consists of one graph $F$, then  
an edge disjoint ${\cal F}$-packing and an edge disjoint ${\cal F}$-factor are called simply an {\em edge disjoint $F$-packing} and an {\em edge disjoint $F$-factor}, respectively. 
Accordingly, the {\em edge disjoint $F$-packing problem} is the problem of finding in $G$ an edge disjoint $F$-packing ${\cal Q}$ having the maximum number of edges $|E({\cal Q})|$ or, equivalently, the maximum number of parts $k({\cal Q})$. Let $\lambda _e(G)$ denote the number $k({\cal Q})$
of parts in a maximum edge disjoint $\Lambda $-packing of $G$.

A graph $D$ is called the {\em line graph of} a simple graph $G$
if $V(D) = E(G)$ and $ab \in E(D)$ if and only if edges $a$ and $b$ in $G$ have a common end-vertex.
Let $L(G)$ denote the line graph of a graph $G$.
A graph $G$ is called a {\em line graph} if there exists a graph $F$ such that $G = L(F)$.
It is known (and easy to show) that if two non-isomorphic graphs $A$ and $B$ are such that $L(A)$  and $L(B)$ are isomorphic, then $\{A, B\} = \{Y, \Delta\}$, where $Y$ is a claw and $\Delta $ is a triangle and both $L(A)$  and $L(B)$ are triangles.
Therefore if $H \not \in \{Y, \Delta\}$ and $H$ is a line graph, then there is a unique graph $F$ such that $H = L(F)$.
A packing $P$ in a graph $G$ is called an 
{\em induced packing} in $G$ if $P$
is an induced subgraph of $G$.
\\[1ex]
\indent
We need the following simple observations.
Let ${\cal F}$ a family of non-isomorphic graphs and 
$L({\cal F}) = \{L(F): F \in {\cal F}\}$.
\bs 
\label{inducedpacking}
Let $G$ be a graph.
If $P$ is an ${\cal F}$-packing in $G$, then $L(P)$ is an induced 
$L({\cal F})$-packing in $L(G)$.
If $Y$ and $\Delta $ are not in ${\cal F}$ and
$L(P)$ is an induced
$L({\cal F})$-packing in $L(G)$, then $P$ is an ${\cal F}$-packing in $G$. 
In particular, if $P$ is an ${\cal F}$-factor in $G$, 
then $L(P)$ is a vertex maximum  induced  
$L({\cal F})$-packing in $L(G)$.
\es

\bs
\label{EdgePackingLineGraph}
Let $G$ be a graph and $D = L(G)$. 
Then the following holds.
\\[0.7ex]
$(a1)$ Let $Q = \{Q_i:  i \in \{1, \ldots , k\}\}$ be an edge 
disjoint ${\cal F}$-packing in $G$. Then 
$L(Q)$ is an ${\cal F}'$-packing in  $D$, where
$Cmp (L(Q)) = \{L(Q_i) : i \in \{1, \ldots , k\}\}$ and
${\cal F}' = \{L(F): F \in  {\cal F}\}$.
\\[0.7ex]
$(a2)$
Let $P$ be a packing in $D$ and 
$Cmp (P) = \{P_i : i \in \{1, \ldots , k\}\}$.
Then $\{Q_i = V(P_i): i \in \{1, \ldots , k\}\}$ is an edge 
disjoint packing in $G$
with $P_i$ being a spanning subgraph of $L(\dot{Q}_i)$, where 
$\dot{Q}_i$ is the subgraph in $G$ induced by edge subset $Q_i$.
\es

These observations allow to deduce various byproducts from the  packing results described before and obtain some facts on edge disjoint packings in a graph. 
Here are some of these results.
\\[1ex]
\indent
From {\bf \ref{inducedpacking}} we have, in particular:
\bs 
\label{inducedmatching}
Let $G$ be a graph.
Then $P$ is a $\Lambda $-packing in $G$ if and only if
$L(P)$ is an induced matching in $L(G)$.
\es

Since the $\Lambda $-packing problem is $NP$-hard even for cubic graphs, we have:
\bs The induced matching problem is $NP$-hard 
for line graphs of cubic graphs.
\es

Using a procedure  for connected graphs similar to Procedure ${\cal E}'$ for 2-connected graphs, it is easy to show the following:
\bs
\label{MatchingInClawfreeGraphs}
Let  $G$ be a connected claw-free graph. Then the following holds.
\\[0.5ex]
$(a1)$ If $M$ is a maximum matching in $G$, then 
$e(M) = \lfloor v(G)/2 \rfloor $.
\\[0.5ex]
$(a2)$ If $v(G) \equiv 1 \bmod 2$, then for every edge $e$ in $D$ there exist a maximum matching $M$ in $G$ that avoids $e$.
\es

Since $L(G)$ is a claw-free graph, we have  from 
{\bf \ref{EdgePackingLineGraph}} and
{\bf \ref{MatchingInClawfreeGraphs}}:
\bs
Let $G$ be a connected graph. Then
\\[0.5ex]
$(a1)$
$\lambda _e(G) = \lfloor e(G)/2 \rfloor $ and
\\[0.5ex]
$(a2)$ if $e(G) \equiv 1 \bmod 2$, then for every 3-vertex path $L$ in $G$ there exists  a maximum edge disjoint 
$\Lambda $-packing ${\cal Q}$ such that $L$ is not a member of ${\cal Q}$.
\es

 An edge disjoint factor ${\cal Q}$ of $G$ is said to be
an {\em edge  $k$-factor} if every member of ${\cal Q}$ induces in $G$  a connected graph having $k$ edges.
\bs
\label{edge3packing}
Suppose that $G$ is a graph such that $L(G)$ is connected and has at most two end-blocks.
If $e(G) \equiv 0 \bmod 3$, then $G$ has an edge 3-factor.
\es

{\bf Proof} (uses {\bf \ref{eb<3}} and 
{\bf \ref{EdgePackingLineGraph}}).
Since $V(L(G)) = E(G)$, we have:
$e(G) \equiv 0 \bmod 3 \Rightarrow v(L(G)) \equiv 0 \bmod 3$.
Since $L(G)$ is claw-free, by {\bf \ref{eb<3}}, $L(G)$ has a 
$\Lambda $-factor. Therefore we are done by 
{\bf \ref{EdgePackingLineGraph}}.
\ep
\\[1ex]
\indent
We call a graph $G$ an {\em edge-chain} if 
$G - Lv(G) = (\cup \{B_i: i \in \{1, \ldots , k\}\}) \cup
 \{e_i:  i \in \{1, \ldots , k-1\}\}$, where each $B_i$ is an edge 2-connected graph, all $B_i$' are disjoint, and each $e_i$ is an edge with one end-vertex in $B_i$ and the other end-vertex in  $B_{i+1}$.
 It is easy to see the following.
 \bs
 \label{edge-chain}
If $G$ is an edge-chain, then $L(G)$ has at most two end-blocks.
 \es

From {\bf \ref{clfree,2con-avoid-e}} $(a1)$,  {\bf \ref{clfree,3con,einL}}, {\bf \ref{edge3packing}}, and {\bf \ref{edge-chain}} we have:
\bs
If $G$ is an edge-chain and  
$e(G) \equiv 0 \bmod 3$,
then $G$ has a 3-edge factor.
Moreover,
\\[0.7ex]
$(a1)$ if $G - Lv(G)$ is edge 2-connected and $e(G) \equiv 0 \bmod 3$, then for every 
3-vertex path $L$ in $G$ there exists an edge 3-factor ${\cal Q}$ with no member containing $L$,
\\[0.7ex]
$(a2)$ if $G$ is edge 3-connected $e(G) \equiv 0 \bmod 3$, then for every 
3-vertex path $L$ in $G$ there exists an edge  3-factor ${\cal Q}$ with a member containing  $L$.
\es

In \cite{KKN} we put forward the following conjecture.
\bs {\em {\bf Conjecture.}}
\label{3conclawfreePi}
Every 3-connected claw-free graph with $v(G) \equiv 0 \bmod 4$ has a $\Pi $-factor.
\es

By {\bf \ref{inducedpacking}}, Conjecture 
{\bf \ref {3conclawfreePi}} is equivalent to the following conjecture on induced $\Lambda $-packings.
\bs {\em {\bf Conjecture.}}
\label{inducedLpacking}
If $G$ is a 3-connected claw-free graph with $v(G) \equiv 0 \bmod 4$ and $P$ is a maximum induced
 $\Lambda $-packing in $L(G)$, then  
 $\lambda (P) = v(G)/4$.
\es



As we mentioned in the introduction, the problem of packing induced 3-vertex paths in a claw-free graph, interesting in itself, is also related to the Hadwiger conjecture.

Let $h(G)$ be the maximum integer $r$ such that $G$ has $K_r$ as a minor.
In 1943 Hadwiger conjectured that if a graph  
$G$ has no proper vertex coloring with $s-1$ colors, then $h(G) \ge s$ (see \cite{D}). 
Now consider a graph $F$  with 
$\alpha (F) = 2$, where $\alpha (F)$ in the size of a maximum vertex subset in $F$ with no two adjacent vertices, and so $F$ is claw-free. Then obviously the vertices of $F$ cannot be colored properly with $s-1$ colors, where $s = \lceil v(G)/2 \rceil$. 
Thus, a natural (open) question  is whether $h(F) \ge s$ as the Hadwiger conjecture claims.
If $P$ is a $\Lambda $-packing in $F$ such that every component 
(3-vertex path) of $P$ is an induced subgraph in $F$, then every two components of $P$ are connected by an edge in $G$. Therefore contracting each component $L$ of $P$ to a new vertex $c(L)$ results in a graph $G'$ having  the complete subgraph $K$ with 
$V(K) = \{c(L): L \in Cmp(P)\}$, and so $v(K) = \lambda (P)$. 
Thus, the maximum packing of induced 3-vertex paths in $F$ provides a maximum complete minor $K$ of $F$ in which 
every vertex corresponds to an induced 3-vertex path in $F$.

Let $h'(F)$ be the maximum integer $r$ such that $F$ has a minor $K_r$ in which every vertex corresponds to either a vertex or an edge in $F$. Obviously, $h(F) \ge h'(F)$.
In 1999 \cite{K-Hadw} we proved that  if $F$ is not $s$-connected, then the Hadwiger conjecture is true, moreover, $h'(F) \ge s$.
\\[1ex]
\indent
From {\bf \ref{EdgePackingLineGraph}} we have in particular:
\bs
\label{inducedLpack} 
Let $G$ be a graph and $P$ a subgraph of $G$.
Then $P$  is an edge disjoint $\Pi$-packing in $G$
if and only if  $L(P)$ is a packing of induced 3-vertex paths in $L(G)$.
\es

It is known \cite{DT} that the maximum edge disjoint $\Pi$-packing problem is $NP$-hard. Therefore by 
{\bf \ref{inducedLpack}}, the maximum packing of induced 3-vertex paths  is also $NP$-hard.
However, probably the following is true.
\bs {\em {\bf Conjecture.}}
\label{alpa=2} Let $\alpha $ be a positive integer.
Then there exists a polynomial-time algorithm $A_{\alpha}$ for finding a maximum  packing of disjoint
induced 3-vertex paths in a claw-free graph $G$ with $\alpha (G) \le \alpha $. 
\es


\end{document}